\documentclass[reqno,11pt]{amsart}
\usepackage[usenames, dvipsnames]{color}
\usepackage{comment}
\usepackage{enumitem}
\usepackage{marginnote}
\usepackage{todonotes}
\usepackage{esint}
\usepackage{graphicx}
\usepackage{amsfonts}
\usepackage{mathrsfs}
\usepackage[english]{babel}
\usepackage{amsmath, amsthm, amssymb, amsbsy}
\usepackage{mathtools}

\usepackage[colorlinks=true]{hyperref}

\newcounter{thmcounter}
\numberwithin{equation}{section}
\numberwithin{thmcounter}{section}
\theoremstyle{definition}

\theoremstyle{plain}
\newtheorem{lemma}[thmcounter]{Lemma}
\newtheorem{thm}[thmcounter]{Theorem}

\newtheorem{prop}[thmcounter]{Proposition}
 
\newtheorem{rmkn}[thmcounter]{Remark}

\theoremstyle{definition}
\newtheorem{definition}[thmcounter]{Definition}

\DeclareMathOperator{\dive}{div}

\newcommand{\wei}[1]{\langle #1 \rangle}
\newcommand{\R}{{\mathbb R}}

\newcommand{\sB}{\mathscr{B}}
\newcommand{\sP}{\mathscr{P}}
\newcommand{\A}{\mathbf A}
\newcommand{\bB}{\mathbf B}
\newcommand{\bI}{\mathbf I}
\newcommand{\dg}{\mathbf{d}g}
\newcommand{\twd}{\widetilde{w}^\delta}

\title[Inviscid limits]{An inviscid limit problem for Navier-Stokes equations in 3D domains with oscillatory boundaries}

\author[T.~Phan]{Tuoc Phan}
\address[T.~Phan]{Department of Mathematics, University of Tennessee, 227 Ayres Hall, 1403 Circle Drive, Knoxville, TN 37996, USA}
\email{tphan2@utk.edu}

\author[D.~Valdebenito]{Dar\'{i}o A.~Valdebenito}
\address[D.~Valdebenito]{Department of Mathematics, Ave Maria University, 5050 Ave Maria Boulevard, Ave Maria, FL 34142, USA}
\email{dario.valdebenito@avemaria.edu}

\thanks{The research of T. Phan is partially supported by Simons Foundation, grant \# 769369.}

\begin{document} 
\begin{abstract} 
We study an inviscid limit problem for a class of Navier-Stokes equations with vanishing measurable viscous coefficients in 3-dimensional spatial domains whose boundaries are oscillatory, depending on a small parameter, and become flat when the parameter converges to zero. Under some sufficient conditions on the anisotropic vanishing rates of the eigenvalues of the matrices of the viscous coefficients and the oscillatory  parameter, we show that Leray-Hopf weak solutions of the Navier-Stokes equations with no slip boundary condition converge to solutions of the Euler equations in the upper half space. To prove the result, we apply a change of variables to flatten the boundaries of the spatial domains for the Navier-Stokes equations, and then construct the boundary layer terms. As the Navier-Stokes equations and the Euler equations are originally written in two different domains, additional boundary layer terms are constructed and their estimates are obtained.
\end{abstract}
\maketitle


{\emph{AMS Classification:}  76D10, 76D05,  35B35, 35Q30, 35Q86, 35B30

\section{Introduction} 
\subsection{Problem Setting and Main Result} We study an inviscid limit problem for inhomogeneous fluids in 3D domains with vanishing viscous coefficients which are matrices of measurable functions. The considered spatial domains have oscillatory boundaries depending on a vanishing parameter $\delta>0$, and they become flat as $\delta \rightarrow 0^+$. To set up the problem, for a given sufficiently small parameter $\delta>0$, let us denote
\begin{equation}\label{eq:Omega}
\Omega=\Omega^\delta= \big \{(x',x_3) \in \R^3: x_3> \delta^\alpha g(x'/\delta), \ x' = (x_1, x_2) \in \R^2 \big\},
\end{equation}
where $g:\R^2 \to [0,\infty)$ is a $C^2$-function, and $\alpha>5/2$ is a fixed constant.  We assume that
\begin{equation}\label{eq:gL}
\sum_{|\gamma|=1}^2 \|D^\gamma g\|_{L^\infty(\R^2)}=L  \in [0, \infty),
 \end{equation}
 where $\gamma=(\gamma_1,\gamma_2) \in (\mathbb{N} \cup \{0\})^2$ is a multi-index, and $|\gamma| = \gamma_1 + \gamma_2$. Most of the time, we omit the superscript $\delta$ on $\Omega$ for the sake of notational simplicity. 

\smallskip
For each $T \in (0, \infty]$, we denote
\[ 
\Omega_T=\Omega_T^\delta = (0, T)\times \Omega. 
\]
Let  $\A_0: \Omega_T \rightarrow \mathbb{R}^{3 \times 3}$ be the matrix of measurable functions representing the viscous coefficients of an inhomogeneous fluid in $\Omega_T$. We assume that $\A_0$ is symmetric and satisfies the following ellipticity and boundedness conditions: there exist $\Lambda \in (0,1)$ and $0 < \nu< \eta < 1$ such that
\begin{equation} \label{ellipticity-cond}
\left\{
\begin{split} 
& \Lambda \big(\eta |\xi'|^2 + \nu |\xi_3|^2\big) \leq\wei{\A_0(t,x)\xi, \xi}, \\
& |\wei{\A_0(t,x)\xi,\zeta}|\leq \Lambda^{-1}\Big(\eta |\xi'| |\zeta'| + \nu \big( |\xi_3| | \zeta'| + |\zeta_3| |\xi'| + |\xi_3| |\zeta_3| \big) \Big),
\end{split}   \right.
\end{equation}
for all $\xi = (\xi', \xi_3) \in \mathbb{R}^2\times \mathbb{R}$ and $\zeta = (\zeta', \zeta_3) \in \mathbb{R}^2\times \mathbb{R}$, and for all $(t, x) \in \Omega_T$.  Note that in \eqref{ellipticity-cond}, and throughout the paper, $\wei{\cdot, \cdot}$ denotes the usual inner product in Euclidean spaces.  It is not difficult to construct a matrix $\A_0$ satisfying \eqref{ellipticity-cond}; for instance, $\A_0=\text{diag}\,(\eta,\eta,\nu)+\mathbf{H}$, where the $L^\infty$ norm of the entries of the symmetric matrix $\mathbf{H}$ is $o(\nu)$ as $\nu\to 0$.

\smallskip
We consider the following system of Navier-Stokes equations of an incompressible inhomogeneous fluid in $\Omega_T$ with a measurable viscosity coefficient matrix $\A_0$, and an unknown fluid velocity $\bar{u}: \Omega_T \rightarrow \R^3$:
\begin{equation}\label{eq:NSE}
\left\{
\begin{array}{cccl}
\partial_t \bar{u} - \dive[\A_0 \nabla \bar{u}]+ (\bar{u} \cdot \nabla) \bar{u} &= &-\nabla \bar{p}+\bar{F}  &  \text{in }\Omega_T,\\
\nabla \cdot \bar{u}&=  & 0 &  \text{in }\Omega_T, \\
\bar{u}&=  &0  &\text{on } (0, T) \times \partial\Omega,\\
\bar{u}(0, \cdot)&=  & \bar{U}(\cdot) & \text{on } \Omega.
\end{array} \right.
\end{equation}
Here, $\bar{U}: \Omega \rightarrow \R^3$ is a given initial data satisfying $\nabla \cdot \bar{U} =0$, $\bar{F}: \Omega_T \rightarrow \R^3$ is a given external force, and $\bar{p}: \Omega_T \rightarrow \R$ is the fluid pressure which is unknown.  It is worth mentioning that the motivation of the study of  \eqref{eq:NSE}  arise from other studies such as geophysical fluids,  oceanic models, thermal fluids, and non-Newtonian fluid models  (see, e.g., \cite{CDGG, FNS, HS,  Lions, O, Ped, Tarfulea}) in which domains are non-flat and the viscosity coefficients are measurable functions.

\smallskip
From the assumption \eqref{ellipticity-cond}, we note that $\A_0$ has eigenvalues of order $\eta$ and $\nu$,  and formally
\begin{equation} \label{A-0-converges-zero}
 \A_0 =\A_0^{\eta, \nu} \rightarrow \mathbf{0} \quad \text{as} \quad (\eta, \nu) \rightarrow (0^+, 0^+).
 \end{equation}
In addition, we observe from the definition of $\Omega^\delta$ in \eqref{eq:Omega} that the spatial domain $\Omega^\delta$ becomes the upper-half flat domain $\Omega^0:=\R_+^3=\R^2\times (0,\infty)$ as $\delta \rightarrow 0^+$. Due to this and \eqref{A-0-converges-zero}, when $(\eta, \nu, \delta) \rightarrow (0^+,0^+,0^+)$, it is a natural question to study if solutions $\bar{u}$ to \eqref{eq:NSE} converge to the solution $w^0$ of the following Euler equations
\begin{equation}\label{eq:Euler}
\left\{
\begin{array}{cccl}
\partial_t w^0+(w^0\cdot\nabla)w^0&= & -\nabla q+F^0 & \text{in }   \Omega_T^0,\\
\nabla\cdot w^0&= & 0 & \text{in }   \Omega_T^0,\\
w^0\cdot \vec{e}_3 &=  &0 &\text{on } (0, T) \times \partial\Omega^0,\\
w^0(0,\cdot)&= &W^0(\cdot)&   \text{on }   \Omega^0.
\end{array} \right.
\end{equation}
Here, in \eqref{eq:Euler}, $\vec{e}_3=(0,0,1)$ is the inward-pointing unit normal vector on $\partial\Omega^0$, $\Omega_T^0=(0,T)\times\R_+^3$, and $q: \Omega_T^0 \rightarrow \R$ is the corresponding unknown pressure. Moreover, $F^0$ and $W^0$ are suitable given forcing term and initial data, respectively.  

\smallskip
The objective of this paper is to establish sufficient conditions on $\A_0$, $\delta$, and $g$ so that such the inviscid limit holds.  To state our main result, let us recall the following classical result on the existence and uniqueness of the strong solution of the Euler equations \eqref{eq:Euler} (see  \cite{Brezis}, \cite{Temam-1}).
\begin{lemma} \label{Euler-existence} 
Let $s>5/2$, $W^0 \in H^s(\mathbb{R}^3_+)^3$, and $F^0 \in L^1((0, T), H^s(\mathbb{R}^3_+)^3)$. Assume that
\[
  \nabla \cdot W^0=0 \quad \text{in} \quad \mathbb{R}^3_+, \quad \text{and} \quad W^0 \cdot \vec{e}_3 =0.
\]
Then, there exist $T^*\in (0,\infty)$ and a unique strong solution 
\[ w^0 \in L^\infty((0, T^*), H^s(\mathbb{R}^3_+)^3)
\]
to \eqref{eq:Euler}. Moreover, $\nabla q\in L^\infty((0,T^*),L^\infty(\mathbb{R}^3_+))$.
\end{lemma}

\smallskip
Next, we provide the set of parameters $(\eta, \nu, \delta)$ required for our result.
\begin{definition}\label{hyp:thm}  For given positive constants $K_0$, $\delta_0$, and for a given continuous function $\beta:[0,1]^2\to (0, \infty)$ satisfying
\begin{equation} \label{beta-fn-def}
\lim_{(\eta,\nu)\to (0^+,0^+)}\beta(\eta,\nu)=0,
\end{equation}
we define 
\[
\begin{split}
E(K_0, \delta_0, \beta) =  \Big \{& (\eta,\nu, \delta)  \in (0,1)^2 \times (0, \delta_0)  :0<\nu<\eta <1,  \\
\quad & \delta^{\alpha-1} \eta \leq K_0 \nu,  \text{ and } \eta+\sqrt{\nu/\eta} \leq  \beta(\eta,\nu)^2 \Big\}.
\end{split}
\]
\end{definition}
\smallskip \noindent
In addition, let us define $\Psi_0:\R_+^3\to\Omega$ as
\[
\Psi_0(y', y_3) = (y', y_3 + \delta^\alpha g(y'/\delta)), \quad (y', y_3) \in \mathbb{R}^2 \times (0, \infty). 
\]
For a given function $G: (0, T) \times \Omega \rightarrow \mathbb{R}^3$, the composed function $G \circ \Psi_0: (0, T) \times \R_+^3 \rightarrow \mathbb{R}^3$ is defined by
\[
G \circ \Psi_0 (t, y) = G(t, \Psi_0(y)), \quad (t, y) \in (0, T) \times \R_+^3.
\] 
Throughout this paper, for each $k=1,2,3$, we denote $D_k u$ the directional derivative of a function $u$ in the $\vec{e}_k$-direction, and $D_{x'} u = (D_1 u, D_2u)$.  We refer readers to Section \ref{def-weak-sol} below for the definition of Leray-Hopf weak solutions  $\bar{u} \in L^\infty((0, T), L^2(\Omega)^3) \cap L^2((0, T), H^1_0(\Omega)^3)$ to \eqref{eq:NSE} in $\Omega_T$ with $T>0$.  
\smallskip
The following theorem is the main result of the paper. 
\begin{thm}\label{thm:main}   Let $\Lambda \in (0,1)$, $\alpha, s \in (5/2, \infty)$, and let $L, K_0$ be positive numbers. Then, there is a sufficiently small number $\delta_0 = \delta_0(\Lambda, L, K_0)>0$ such that the following assertion holds. Assume that \eqref{eq:gL} and \eqref{ellipticity-cond} are satisfied with $(\eta, \nu, \delta) \in E(K_0, \delta_0, \beta)$ for some function $\beta$ satisfying \eqref{beta-fn-def}. Let  $F^0\in L^1((0,T), H^s(\R_+^3)^3)$, and let $W^0 \in H^{s}(\R_+^3)^3$  satisfy
\[
  \nabla \cdot W^0 =0 \quad \text{in } \R_+^3 \quad \text{and} \quad W^0\cdot \vec{e}_3 =0 \quad \text{on } \partial\R_+^3.
\]
Also, let $\bar{U}  \in L^2(\Omega)^3$ and $\bar{F} \in  L^2((0,T^*),L^2(\Omega)^3)$ satisfy $\nabla \cdot \bar{U}  =0$ in $\Omega$ and
\begin{equation} \label{ini-cond}
\begin{split}
& \|\bar{U}\circ \Psi_0 -W^0\|_{L^2(\R_+^3)} \leq \beta +\delta^{\alpha-5/2},\\
& \|\bar{F}\circ \Psi_0 -F^0\|_{L^1((0,T^*),L^2(\R_+^3))} \leq \beta + \delta^{\alpha-5/2},
\end{split}
\end{equation}
where $T^* \in (0, T]$ is defined in  \textup{Lemma \ref{Euler-existence}}. Then, for every Leray-Hopf weak solution
\[
\bar{u} \in L^\infty((0, T^*), L^2(\Omega)^3) \cap L^2((0, T^*), H^1_0(\Omega)^3)
\]
of \eqref{eq:NSE} in $\Omega_{T^*}$, it holds that
\begin{equation}\label{eq:uwbound1}
\|\bar{u} \circ \Psi_0 -w^0\|_{L^\infty((0,T^*),L^2(\R_+^3))}\leq C(\beta+\delta^{\alpha-5/2}),
\end{equation}
and
\begin{equation}\label{eq:uwbound2}
\eta \| D_{x'} \bar{u}\|_{L^2(\Omega_{T^*})} + \nu \| D_{3} \bar{u}\|_{L^2(\Omega_{T^*})}  \leq C (\beta+\delta^{\alpha-5/2}).
\end{equation}
Here, in \eqref{eq:uwbound1} and \eqref{eq:uwbound2}, $w^0 \in L^\infty((0, T^*), H^s(\R_+^3)^3)$ is the unique strong solution to \eqref{eq:Euler} whose existence is assured in \textup{Lemma \ref{Euler-existence}},  and
\[ C= C(T^*, \|w^0\|_{L^\infty((0,T^*),H^s(\R_+^3))}, \|F^0\|_{L^1((0,T^*),H^s(\R_+^3))}, L, \Lambda, K_0, s) >0\]
is a constant independent of $\bar{u}$, $\eta$, $\nu$, $\delta$, and $\A_0$.  
\end{thm}

Theorem \ref{thm:main} proves that  even with just measurable viscosity coefficients, a suitable anisotropic convergence of the viscosity, compounded with a suitable convergence of the domain $\Omega^\delta$ to the upper half space as $\delta \rightarrow 0^+$, provides the mechanism for the convergence of the velocity $\bar{u}$ of a viscous fluid in $\Omega_{T^*}^\delta$ to the velocity $w^0$ of an ideal fluid in the upper half space (or, more precisely, the convergence of $\bar{u}\circ\Psi_0$ to $w^0$).  The novelty in  Theorem \ref{thm:main} is that it does not require any regularity conditions on the matrix $\A_0$, rather, we only require the uniform ellipticity and boundedness conditions \eqref{ellipticity-cond}.  As such, our result could be useful in other problems such as \cite{Liu, Tarfulea} in which an incompressible fluid is coupled with some other physical quantities, and in the studies of  geophysical fluids  such as in  \cite{CDGG, FNS, O, Ped}, and in other studies such as \cite{AGZ, DP, DPK, HS, Hou-Huang, Lions} in which the viscosity coefficients are measurable functions.

\smallskip
We remark that when $\delta=0$ and $g\equiv 0$, Theorem \ref{thm:main} recovers results in \cite{Mas}, where $\A_0=\text{diag}(\eta,\eta,\nu)$ is a constant, $3\times 3$ diagonal matrix.  It is also worth pointing out that a similar problem, but in 2D and with $\A_0 = \text{diag}(\nu, \nu)$ (diagonal $2\times 2$ matrix) and the Navier slip boundary condition, was studied in \cite{GLNR}: under suitable conditions on the oscillating parameter and the viscosity, the authors showed that solutions of the Navier-Stokes system converge to solutions of the Euler system. Besides working on 3D domains, our approach to prove Theorem \ref{thm:main} is to flatten the domain, and then construct suitable boundary layers to derive the energy estimates, which is different from the approach used in \cite{GLNR}, which relies on stability analysis. 

\smallskip
We also remark that an inviscid limit problem for Navier-Stokes equations with measurable viscosity coefficient matrix $\A_0$ that satisfies some slightly different conditions compared to \eqref{ellipticity-cond}, and  on domains with sufficiently smooth boundaries is studied in \cite{PV}. Under suitable anisotropic conditions on the convergence $\A_0 \rightarrow \mathbf{0}$, weak solutions of the Navier-Stokes equations are proved to converge to solutions of the Euler equations. Compared to \cite{PV}, in this paper our  considered domains are not fixed, as they depend on the parameter $\delta$ which also convergences to zero. In addition, unlike \cite{PV} in which the Navier-Stokes equations and Euler equations are in the same domains, Theorem \ref{thm:main} provide a convergence result in the case that the two equations are posed in different domains.

\subsection{Some other related literature} The literature on the \emph{inviscid limit} problem is vast and here we only briefly mention a few related works. In the case of the entire space and  with constant diagonal viscous coefficients, the problem was solved in the well known papers \cite{Kato, Swan}. In the case where the domains have non-empty flat boundaries, the initiation of the study was given in \cite{Kato-1} in which sufficient conditions on the boundary layers are given. Some more recent examples can be found in \cite{Mas, Wang} where the boundaries of the domains considered are also flat, and in \cite{DN} for a recent interesting work in bounded domains. In addition, interested readers may find in \cite{D-GV, SV, KNVW, Li-Wang, Liu-Yang, Maekawa, SC-1, SC-2}, the survey paper \cite{MM}, and references therein, both recent and classical results on the well-posedness and regularity of solutions of the Prandtl equations in upper-half spaces. Some other results on inviscid limit problems can be found, for instance, in \cite{GHT, Wang} in which the viscosity coefficients are assumed to be a constant diagonal matrix.

\subsection{Main ideas and approaches}
The proof of Theorem \ref{thm:main} begins with flattening the domain $\Omega$ by taking $u = \bar{u} \circ \Psi_0$, so the system of Navier-Stokes equations \eqref{eq:NSE} becomes the system \eqref{NS.eqn} of $u$ in the upper-half space $\R_+^3=\R^2\times (0,\infty)$ with measurable coefficients. We construct explicitly a corrector $\sB$, or \emph{boundary layer}, adjusting ideas from \cite{Mas} and \cite{PV}; in particular, the corrector involves a perturbation $w$ of the strong solution $w^0$ of the Euler equations, rather than $w^0$ itself: see the beginning of Section \ref{sect-3.1} for more details and intuition, and \eqref{w-short-def-12-31} for the definition of $w$. The system of equations of $u- w - \sB$ has coefficients depending on the variables $(t,y)\in (0,T)\times\R_+^3$ and on the coefficients $\eta$, $\nu$, and $\delta$, so energy analysis estimates are developed in Sections \ref{B-sec}--\ref{Ener-sec} to control the error terms arising from the change of variables, the perturbation $w$, and the specifics of the setting in the present paper. We then close the energy estimate by applying a Gr\"{o}nwall type lemma.

\subsection{Organization of the paper}
The rest of the paper is organized as follows. In Section \ref{flat-n}, we provide the analysis of flattening the domain and some preliminary analysis needed in the rest of the paper. In particular, the second order operator arising from the change of variables to flatten the domain is studied in Subsection \ref{flat-domain}. In Section \ref{B-sec}, we provide the  construction of the corrector term $w$ for $w^0$, the boundary layer $\sB$, and their estimates. The energy estimates of $u- w - \sB$ are given in Section \ref{Ener-sec}. The proof of Theorem \ref{thm:main} is given in Section \ref{proof-thrm} using ingredients developed in Sections \ref{flat-n}--\ref{Ener-sec}.

\section{Flattening the domains and preliminary analysis} \label{flat-n}
\subsection{Flattening the domains} \label{flat-domain} Note that equations \eqref{eq:NSE} and \eqref{eq:Euler} are posed in two different domains. To compare the two systems, we transform the equations \eqref{eq:NSE} into equations on the upper-half space in which \eqref{eq:Euler} is written. We  also provide some important estimates on the second order operator of the resulting equations. 

\smallskip
Recall that $\Omega=\Omega^\delta$ is the domain defined in \eqref{eq:Omega}, where $g:\R^2\to [0,\infty)$ satisfies \eqref{eq:gL}. We denote $\R^3_+=\R^2\times (0,\infty)$, and let $\Phi_0: \Omega \rightarrow \R^3_+$ be defined by
\[
y = \Phi_0 (x) = (x',x_3-\delta^\alpha g(x'/\delta)), \quad x = (x', x_3) \in \Omega.
\]
Also, let $\Psi_0 = \Phi_0^{-1} : \R^3_+ \rightarrow \Omega$, that is,
\begin{equation}\label{eq:Psi0}
x= \Psi_0(y)= (y', y_3+ \delta^{\alpha} g(y'/\delta)), \quad y = (y', y_3) \in \R^3_+,
\end{equation}
the same function considered in the paragraph following Definition \ref{hyp:thm}. A simple calculation shows that
\[
 D\Phi_0(x)  =  \begin{bmatrix}
1 & 0 & 0\\
0 & 1 & 0\\
-\delta^{\alpha-1}D_1g(x'/\delta) & -\delta^{\alpha-1}D_2g(x'/\delta) & 1
\end{bmatrix},
\]
and
\[
D\Psi_0(y)  = \begin{bmatrix}
1 & 0 & 0\\
0 & 1 & 0\\
\delta^{\alpha-1}D_1g(y'/\delta) & \delta^{\alpha-1}D_2g(y'/\delta) & 1
\end{bmatrix}.
\]
Here, $D_k$ denotes the directional derivative along the standard Euclidean unit vector $\vec{e}_k$ for $k =1,2$. Note that
\[
D\Phi_0 = [D\Psi_0]^{-1} \quad \text{and} \quad \text{det} (D\Phi_0) = \text{det} (D\Psi_0) =1.
\]

\smallskip
Let us also denote
\begin{equation*}
\bB = \bB(x') = \begin{bmatrix}
1 & 0 & 0\\
0 & 1 & 0\\
-\delta^{\alpha-1}D_1g(x'/\delta) & -\delta^{\alpha-1}D_2g(x'/\delta) & 1
\end{bmatrix}, \quad x ' \in \R^2.
\end{equation*}
We note that $\bB(x') =  D\Phi_0(x)$. Also, by \eqref{eq:Psi0}, it follows that $D\Phi_0(y) = D\Phi_0(\Psi_0(y))$ for all $y = (y', y_3) \in \mathbb{R}^3$. Therefore, in many places, we also write $\bB = \bB(y) =D\Phi_0(\Psi_0(y)) =\bB(y')$ with
\begin{equation}\label{eq:bBy}
\bB(y')  =
\begin{bmatrix}
1 & 0 & 0\\
0 & 1 & 0\\
-\delta^{\alpha-1}D_1g(y'/\delta) & -\delta^{\alpha-1}D_2g(y'/\delta) & 1
\end{bmatrix}, \quad y' \in \mathbb{R}^2,
\end{equation}
as $\bB$ does not depend on $y_3$---but we will often omit the dependence of $\bB$ on $y'$ or $y$ for notational simplicity.  Also, the matrix $\bB$ can be decomposed as 
\begin{equation}\label{eq:bBexp}
\bB=\bI+\delta^{\alpha-1}\bB_g,
\end{equation} 
where $\bI$ the $3\times 3$ identity matrix, and
\begin{equation}\label{eq:bBg}
\bB_g=\begin{bmatrix}
0 & 0 & 0\\
0 & 0 & 0\\
-D_1g(y'/\delta) & -D_2g(y'/\delta) & 0
\end{bmatrix}.
\end{equation}
It will be useful to note that, owing to the fact that $\alpha>1$ and $\delta \in (0,1)$,
\begin{equation}\label{eq:bBbound}
 \|\bB\|_{L^\infty(\R^2)}\leq \max\{L, 1\},
\end{equation}
where $L$ is as in \eqref{eq:gL}. 

\smallskip
Next, let us denote
\begin{equation}\label{eq:bA}
\A(t,y) = \bB(y) \A_0(t, \Psi_0(y)) \bB^*(y) , \quad (t, y) \in (0, T) \times \mathbb{R}^3_+,
\end{equation}
where $ \bB^*$ denotes the adjoint matrix of the matrix $\bB$. As with $\bB$, we may omit the dependence of $\A$ and $\A_0$ on $t$ and $y$ for the sake of simplicity. Using the change of variables 
\begin{equation} \label{change-variables}
\begin{aligned}
u(t, y) &= \bar{u}(t, \Psi_0(y)),& F(t,y) &=\bar{F}(t, \Psi_0(y)), \\
U(y) &= \bar{U}(\Psi_0(y)), & p(t,y) &= \bar{p}(t, \Psi_0(y)),  \quad (t, y) \in (0, T) \times \mathbb{R}^3_+,
\end{aligned}
\end{equation}
and letting
\[
\mathcal{L}(u) =\partial_t u - \dive[\A(t,y) \nabla u] + [(\bB(y) u) \cdot \nabla ] u,
\]
the system of Navier-Stokes equations \eqref{eq:NSE} is transformed to
\begin{equation} \label{NS.eqn}
\left\{
\begin{array}{cccl}
\mathcal{L}(u) + \bB^* \nabla p & =  & F & \text{in }  (0, T) \times \R^3_+,   \\
\dive[\bB u] & = & 0 &\text{in }   (0, T) \times \R^3_+, \\
u & = & 0 &   \text{on }  (0, T) \times \partial \R^3_+, \\
 u(0, \cdot) & =  & U &\text{on } \R^3_+.
\end{array}
\right.
\end{equation}
 
\smallskip
It is simple to prove that the standard energy estimate holds for \eqref{NS.eqn}. Indeed, we formally have
\[
\int_{\R^3_+}\wei{\bB^* \nabla  p, u} dy = \int_{\R^3_+}\wei{\nabla  p, \bB u} dy = - \int_{\R^3_+}  p \dive(\bB u) dy =0.
\]
Here $\wei{\cdot, \cdot}$ denotes the usual inner product in $\R^3$. The same notation will also be used later on for the inner product in $\R^2$. Similarly,
\begin{align*}
\int_{\R^3_+}\wei{[(\bB u) \cdot \nabla ] u , u} dy & =\sum_{k=1}^3 \int_{\R^3_+} \wei{(\bB u)_k D_k u , u}\ dy \\
& = \sum_{k=1}^3\frac{1}{2}\int_{\R^3_+} (\bB u)_k D_k (|u|^2)\ dy\\
& = \frac{1}{2}\int_{\R^3_+} \wei{\bB u, D (|u|^2)}\ dy\\
& = -\frac{1}{2}\int_{\R^3_+}|u|^2 \dive(\bB u)\ dy\\
& =0.
\end{align*}
From the computations above, one can derive the energy estimate for the system  \eqref{NS.eqn}, provided the uniform ellipticity and boundedness of the matrix $\A(t,y)$ have been established. In the proposition below, we state and prove that, under \eqref{eq:gL}, \eqref{ellipticity-cond}, and when $\delta$ is sufficiently small, the viscosity coefficient matrix $\A(t,y)$ appearing in \eqref{NS.eqn} is uniformly elliptic and bounded with eigenvalues of order $\eta$ and $\nu$. 

\begin{prop}\label{prop:Z} 
Let $L$, $K_0$, $\Lambda$ be positive numbers. Then, there exist $\delta_0 = \delta_0 (\Lambda, L, K_0)>0$ sufficiently small and $C = C(\Lambda, L, K_0)>0$ such that if \eqref{eq:gL} and \eqref{ellipticity-cond} hold, and $\delta \in (0, \delta_0)$ satisfies $\delta^{\alpha-1} \eta \nu^{-1} \leq K_0$, then the following statements hold:
\begin{itemize}
\item[\textup{(a)}] For any $\xi=(\xi',\xi_3)\in \R^2\times\R=\R^3$,
\begin{equation}\label{eq:e1}
\frac{\Lambda}{2}\left(\eta|\xi'|^2+\nu|\xi_3|^2\right)\leq \wei{\A_0\bB^*\xi,\bB^*\xi}
\end{equation}
and
\begin{equation}\label{eq:e2pp}
\wei{\A_0\bB^*\xi,\bB^*\xi}\leq C\big(\eta|\xi'|^2+\nu|\xi_3|^2\big).
\end{equation}
\item[\textup{(b)}] For any $\xi=(\xi',\xi_3)$ and $\zeta=(\zeta',\zeta_3)$,
\begin{align*}\label{eq:e2p}
& |\wei{\A_0\bB^*\xi,\bB^*\zeta}|\\ 
& \leq C\left(\eta|\xi'||\zeta'|+\nu|\xi_3\zeta_3|+(\delta^{\alpha-1}\eta +\nu) [|\xi_3||\zeta'|+|\zeta_3||\xi'|]\right).
\end{align*}
\end{itemize}
\end{prop}
\begin{proof}
We begin by proving (b). Recall that from \eqref{eq:bBexp}, we have
\begin{equation*}
\bB=\bI+\delta^{\alpha-1}\bB_g.
\end{equation*} 
where $\bB_g$ is defined in \eqref{eq:bBg}. It then follows that
\begin{align} \notag
\wei{\A_0\bB^*\xi,\bB^*\zeta} & = \wei{\A_0 \xi, \zeta} + \delta^{\alpha-1}\big( \wei{\A_0 \bB_g^* \xi, \zeta} + \wei{\A_0 \xi, \bB_g^* \zeta} \big)\\ \label{0923-24-0}
& \qquad + \delta^{2\alpha-2} \wei{\A_0 \bB_g^* \xi, \bB_g^* \zeta}.
\end{align}
Let us denote $\dg=(-D_1g,-D_2g,0)=(-Dg,0)$. Then, we see that
\[
\bB_g^* \xi = \xi_3 \dg \quad \text{and} \quad \bB_g^* \zeta = \zeta_3 \dg.
\]
From this, and as $\A_0$ is symmetric, we obtain
\begin{align*}
 \wei{\A_0 \bB_g^* \xi, \zeta} + \wei{\A_0 \xi, \bB_g^* \zeta} &= \xi_3 \wei{\A_0 \dg, \zeta} + \zeta_3\wei{\A_0 \dg, \xi}, \\
 \wei{\A_0 \bB_g^* \xi, \bB_g^* \zeta} &= \xi_3  \zeta_3 \wei{\A_0 \dg, \dg}.
\end{align*}
Then, it follows from \eqref{0923-24-0} that 
\begin{align} \notag
\wei{\A_0\bB^*\xi,\bB^*\zeta} &= \wei{\A_0 \xi, \zeta} + \delta^{\alpha-1}\big( \xi_3 \wei{\A_0 \dg, \zeta} + \zeta_3\wei{\A_0 \dg, \xi} \big) \\ \label{0923-24-1}
& \qquad  +  \delta^{2\alpha-2}\xi_3  \zeta_3 \wei{\A_0 \dg, \dg}.
\end{align}
Next, using the second inequality in \eqref{ellipticity-cond},  and \eqref{eq:gL}, 
\begin{align*}
& |\xi_3 \wei{\A_0 \dg, \zeta} + \zeta_3\wei{\A_0 \dg, \xi}| \leq  |\xi_3|| \wei{\A_0 \dg, \zeta}| +| \zeta_3||\wei{\A_0 \dg, \xi}|\\
& \leq \Lambda^{-1}  |\xi_3| \big(\eta |Dg| |\zeta'| + \nu |\zeta_3| |Dg|  \big) + \Lambda^{-1} |\zeta_3|\big( \eta |\xi'| |Dg| + \nu |\xi_3||Dg| \big) \\
& \leq \Lambda^{-1} L \Big( \eta \big(|\zeta_3||\xi'| + |\xi_3||\zeta'|\big) + 2\nu  |\xi_3||\zeta_3|   \Big).
\end{align*}
Similarly, we also have
\begin{align*}
|\xi_3  \zeta_3 \wei{\A_0 \dg, \dg}|  \leq \Lambda^{-1} L^2 \eta |\xi_3| |\zeta_3|.
\end{align*}
Substituting the two foregoing inequalities in \eqref{0923-24-1} and using the second inequality in \eqref{ellipticity-cond}, we obtain
\begin{align*}
\Lambda |\wei{\A_0\bB^*\xi,\bB^*\zeta}| & \leq \Big(\eta |\xi'| |\zeta'| + \nu \big( |\xi_3| \zeta'| + |\zeta_3| |\xi'| + |\xi_3| |\zeta_3| \big) \Big) \\
& \qquad +  L \delta^{\alpha-1}\Big( \eta \big(|\zeta_3||\xi'| + |\xi_3||\zeta'|\big) + 2\nu  |\xi_3||\zeta_3|   \Big)  \\
& \qquad + L^2 \delta^{2\alpha-2}\eta |\xi_3| |\zeta_3|\\
& \leq  \Big(\eta |\xi'| |\zeta'| + \nu\big( 1 + 2L \delta^{\alpha-1}  + L^2 \delta^{2\alpha-2} \eta \nu^{-1}\big)|\xi_3| |\zeta_3| \Big) \\
& \qquad + \big(\nu+  L \delta^{\alpha-1} \eta \big) \big(|\xi_3| \zeta'| + |\zeta_3| |\xi'| \big).
\end{align*}
As $\delta^{\alpha-1} \nu^{-1} \eta \leq K_0$, we can find a constant $C = C(\Lambda, L, K_0)>0$ such that
\begin{align*}
&  |\wei{\A_0\bB^*\xi,\bB^*\zeta}|  \\
 &  \leq C \Big[ \eta |\xi'| |\zeta'|  + \nu |\xi_3||\zeta_3|+ (\delta^{\alpha-1}\eta +\nu) [|\xi_3||\zeta'|+|\zeta_3||\xi'|] \Big].
\end{align*}
Hence (b) is proved.

\smallskip
We now prove \eqref{eq:e2pp}.  Using (b) and the fact that $\delta^{\alpha-1} \eta \nu^{-1} \leq K_0$ yields
\[ 
  |\wei{\A_0\bB^*\xi,\bB^*\xi}|    \leq C \Big[ \eta |\xi'|^2  + \nu |\xi_3|^2+ \nu |\xi_3||\xi'| \Big].
\]
Then, by using the Cauchy-Schwarz inequality and the fact that $\nu < \eta$, we obtain
\[
 |\wei{\A_0\bB^*\xi,\bB^*\xi}|    \leq C \Big[ \eta |\xi'|^2  + \nu |\xi_3|^2 \Big],
\]
where $C = C(\Lambda, L, K_0)>0$. Hence, the second assertion in (a) is proved.

\smallskip
It remains to prove \eqref{eq:e1}. Using \eqref{0923-24-1} with $\zeta=\xi$, we see that 
\begin{equation}  \label{0923-24-2}
\begin{split}
\wei{\A_0\bB^*\xi,\bB^*\xi}  & = \wei{\A_0 \xi, \xi} + 2\delta^{\alpha-1} \xi_3 \wei{\A_0 \dg, \xi}   \\
& \qquad  +  \delta^{2\alpha-2}\xi_3 ^2 \wei{\A_0 \dg, \dg}.
\end{split}
\end{equation}
From this,  and the assumptions  \eqref{eq:gL} and \eqref{ellipticity-cond}, we see that
\begin{align*}
\wei{\A_0\bB^*\xi,\bB^*\xi}   & \geq \Lambda(\eta|\xi'|^2 + \nu |\xi_3|^2) -2 \delta^{\alpha-1} \Lambda^{-1} L |\xi_3|  \big(\eta |\xi'| + \nu |\xi_3| \big)\\
& \geq  \Lambda\Big (\eta|\xi'|^2 + \nu \big(1-2\delta^{\alpha-1} \Lambda^{-2} L \big)|\xi_3|^2 \Big) -2 \delta^{\alpha-1} \eta \Lambda^{-1} L |\xi_3|   |\xi'| 
\end{align*}
where we used the fact that $\wei{\A_0 \dg, \dg} \geq 0$, and that the third component of $\dg$ is zero. Note that it follows from  Cauchy-Schwarz inequality that
\[
2 \delta^{\alpha-1} \Lambda^{-1} L \eta  |\xi_3|   |\xi'|  \leq  \frac{1}{2} \Lambda \eta |\xi'|^2 + 2 \eta \delta^{2\alpha-2} \Lambda^{-4}L^2 |\xi_3|^2.
\]
Then, 
\begin{align*}
\wei{\A_0\bB^*\xi,\bB^*\xi} &  \geq \frac{\Lambda}{2} \Big( \eta |\xi'|^2 +  2\nu \big(1-2\delta^{\alpha-1}\big[ \Lambda^{-2}  + \eta \nu^{-1} \delta^{\alpha-1} \Lambda^{-5}L^2\big] \big)|\xi_3|^2 \Big) \\
& \geq \frac{\Lambda}{2} \Big( \eta |\xi'|^2 +  2\nu \big(1-2\delta^{\alpha-1}\big[ \Lambda^{-2}  + K_0 \Lambda^{-5}L^2\big] \big)|\xi_3|^2 \Big).
\end{align*}
We can choose $\delta_0 = \delta_0(\Lambda, L, K_0)>0$ sufficiently small so that
\[
\delta^{\alpha-1}\big[ \Lambda^{-2}  + K_0 \Lambda^{-5}L^2\big] \leq \frac{1}{4}
\]
holds for all $\delta \in (0, \delta_0)$, which implies
\[
\wei{\A_0\bB^*\xi,\bB^*\xi} \geq  \frac{\Lambda}{2}\Big( \eta|\xi'|^2 + \nu |\xi_3|^2\Big).
\]
Hence \eqref{eq:e1} is proved, and the proof of the proposition is thus completed.
\end{proof}
\subsection{ Leray-Hopf weak solutions}  \label{def-weak-sol}

Let us recall the definition of Leray-Hopf weak solutions for \eqref{eq:NSE}. Let $\A_0: \Omega_T \rightarrow \mathbb{R}^{3\times 3}$ be a $3\times 3$ measurable matrix that satisfies the condition \eqref{ellipticity-cond}, $T \in (0, \infty]$, $\bar{U}  \in L^2(\Omega)^3$ with $\nabla \cdot \bar{U}  =0$ in the weak sense, and $\bar{F} \in L^2((0,T),L^2(\Omega)^3)$. A function $\bar{u} \in L^\infty((0,T),  L^2(\Omega)^3) \cap L^2((0,T), H^1_0(\Omega)^3)$  is said to be a Leray-Hopf weak solution of \eqref{eq:NSE} if the following statements hold:
\begin{itemize}
\item[\textup{(i)}] For every $\phi \in C^\infty_c(\Omega_T)^3$ (the space of smooth and compactly supported functions) such that $\nabla \cdot \phi =0$, we have
\[
\int_{\Omega_T}\big[-\wei{\bar{u},\phi_t} + \wei{\A_0 \nabla \bar{u}, \nabla \phi}\big]\ dxdt = \int_{\Omega_T} \wei{\bar{F}(t,x),\phi(t,x)}\ dx dt.
\]
\item[\textup{(ii)}] For every $\tau \in (0,T)$, there is $C = C(\tau, \nu, \eta, \Lambda) >0$ such that
\[
\begin{split}
& \sup_{t \in (0,\tau)}\int_{\Omega} |\bar{u}(t,x)|^2\ dx + \int_{\Omega_\tau} |\nabla \bar{u}(t,x)|^2\ dxdt \\
& \leq C\left[ \int_{\Omega_\tau} |\bar{F}(t,x)|^2\ dxdt + \int_{\Omega} |\bar{U}(x)|^2\ dx\right].
\end{split}
\]
\item[\textup{(iii)}] For every $\phi \in L^2(\Omega)^3$, the map $t \mapsto \int_{\Omega} \wei{\bar{u}(t,x),\phi(x)}\ dx$ is continuous on $[0, T)$, and for every $\varphi \in C^\infty_c(\Omega)$
\[
\int_{\Omega} \wei{\bar{u}(t,x),\nabla \varphi(x)}\ dx =0, \quad \text{for } t \in (0, T).
\]
\item[\textup{(iv)}] The initial condition is satisfied in the following sense:
\[
\lim_{t\rightarrow 0^+} \|\bar{u}(t, \cdot) - \bar{U}(\cdot)\|_{L^2(\Omega)} =0.
\]
\end{itemize}
We note for each fixed $\eta$, $\nu \in (0,1)$ and $\delta$, $\Lambda>0$, under assumption \eqref{ellipticity-cond}, the matrix $\A_0$ is bounded and uniformly elliptic. Therefore, the existence of a Leray-Hopf weak solution $\bar{u}$ of \eqref{eq:NSE} can be proved using the standard Galerkin method, see \cite{Temam-2} for instance.  Using the change of variables \eqref{eq:Psi0}, we obtain a weak solution $u$ of \eqref{NS.eqn}.
 
\subsection{Helmholtz-Leray projections} We now recall the definition of the Helmholtz-Leray projection  $\sP$ in the upper-half space $\mathbb{R}^3_+$. For $\sigma \in (0, \infty)$, denoting $H^\sigma(\R_+^3)$ the usual Bessel potential spaces, we define
\begin{equation*}
V^\sigma=\{ u\in H^\sigma(\R_+^3)^3:  \dive u=0\text{ and } u_3(y',0)=0\},
\end{equation*}
where $ \dive u=0$ is understood in the weak sense, i.e.,
\[
\int_{\R^3_+} \wei{u(y), \nabla \varphi(y)}\, dy =0, \quad \forall \varphi \in C_c^\infty(\R^3_+),
\]
and $u_3(y',0)=0$ is understood in the sense of trace. We also write
\[
V^0 =\{u\in L^2(\R_+^3)^3:  \dive u=0\},
\]
in which $\dive u=0$ is understood in the weak sense. 

\smallskip
We denote $\sP$ the Helmholtz-Leray projection of $L^2(\R_+^3)^3$ onto $V^0$. To define $\sP$, for each given vector field $f \in H^\sigma(\R_+^3)^3$ with some $\sigma \geq 0$, we seek a function $q_0 : \R^3_+ \rightarrow \R$ solving
\begin{equation} \label{dec-sys}
\left\{
\begin{aligned}
\dive(\nabla q_0 - f) & = 0 & \quad \text{in }\R^3_+, \\
\wei{\nabla q_0 - f), \vec{e}_3} &= 0  & \quad \text{on }\partial\R^3_+,
\end{aligned} \right.
\end{equation}
and satisfying the estimate
\begin{equation} \label{q-est}
 \|\nabla q_0 \|_{H^\sigma(\R_+^3)} \leq C \|f\|_{H^\sigma(\R_+^3)},
\end{equation}
where $\vec{e}_3 =(0,0,1)$. Note that the existence of $q_0$ satisfying \eqref{dec-sys}--\eqref{q-est} with $\sigma=0$ follows from the standard method using the Lax-Milgram theorem. By the classical regularity theory for elliptic equations (see \cite[Chapter 6]{Evans}, for instance), we conclude that \eqref{q-est} also holds for any $\sigma>0$. 

\smallskip
Now, we define  $\sP: H^\sigma(\R_+^3) \rightarrow V^\sigma$ as 
\[ \sP(f) = f -\nabla q_0,  \quad f \in H^\sigma(\R_+^3).
\] 
From \eqref{q-est} and the triangle inequality we obtain, for any $\sigma\geq 0$,
\begin{equation} \label{Projection}
\|\sP(f)\|_{H^\sigma(\R_+^3)}\leq C\|f\|_{H^\sigma(\R_+^3)},
\end{equation}
where $C>0$ is a constant depending only on $\sigma$.

\section{Boundary layers and their estimates}  \label{B-sec}

 Let $u$ be a weak solution of the Navier-Stokes equations \eqref{NS.eqn} and $w^0$ be the strong solution of the Euler equations \eqref{eq:Euler}. Owing to the respective boundary conditions in \eqref{NS.eqn} and \eqref{eq:Euler}, in general we have
\begin{equation} \label{diff-bdr-cond}
(u-w^0)\big|_{\{y_3=0\}} \neq 0.
\end{equation}
Hence, we cannot eliminate the boundary terms in the integration by parts when working on energy estimates for the equations for $u-w^0$. To deal with this, we introduce a corrector 
\begin{equation} \label{hat-sB-term}
\hat{\sB}:=\delta^{\alpha -5/2} \twd + \sB,
\end{equation}
where $\twd$ is defined in \eqref{eq:wdeltadef} and $\sB$ is defined in \eqref{BL.def}.  The rest of the section is devoted to giving the definitions of $\twd$  and $\sB$, and proving their essential estimates.
\subsection{Construction of $\twd$ and its estimates} \label{sect-3.1} Observe that originally, $\bar{u}$ and $w^0$ are defined in different domains, and we use the change of variables to make the domain of $\bar{u}$ the same as that of $w^0$. To construct a corrector, we find an adjustment from $w^0$ to seal the discrepancy between the two domains. For sufficiently small $\delta>0$, and for the given solution $w^0$ of \eqref{eq:Euler} defined on $\Omega_{T^*}^0 = (0,T^*) \times \mathbb{R}^3_+$ with some $T^*>0$, let 
\[ \twd = (\twd_1, \twd_2, \twd_3): \Omega_{T^*}^0 \rightarrow \mathbb{R}^3
\] 
be defined as
\begin{equation}\label{eq:wdeltadef}
\begin{aligned}
\twd_1(t, y)&=-w_1^0(t, y), \quad \twd_2(t, y)=-w_2^0(t, y), \quad \text{and} \\
\twd_3(t, y)&=-w_3^0(t, y)\\
& +(\delta^{3/2}+\delta^{\alpha-1})\left(-D_1g(y'/\delta)w_1^0(t, y)-D_2g(y'/\delta)w_2^0(t, y)\right).
\end{aligned}
\end{equation}
Also, let $w= w^\delta: \Omega_T^0 \rightarrow \mathbb{R}^3$ be the perturbation of $w^0$ which is defined as
\begin{equation}\label{eq:w-def}
w = w^0-\delta^{\alpha-5/2}\twd.
\end{equation}
We can write $\twd$ and $w$ using \eqref{eq:bBexp}:
\begin{equation} \label{w-short-def-12-31}
\begin{split}
& \twd =\big( (\delta^{3/2} + \delta^{\alpha-1}) \bB_g - \bI \big)w^0 \quad \text{and} \\
& w   = \big[\bI + \delta^{\alpha-5/2}\big(\bI - (\delta^{3/2} + \delta^{\alpha-1}) \bB_g \big) \big] w^0.
\end{split}
\end{equation}
Note that we omit the superscript $\delta$ in $w$ for the sake of notational simplicity. In addition, when $\delta=0$, the definition in $\eqref{eq:w-def}$ recovers the solution $w^0$ of \eqref{eq:Euler}.

\smallskip
 In the next lemma we state and prove some equalities and estimates concerning the functions $\twd$ and $w$ defined in \eqref{eq:wdeltadef}, \eqref{eq:w-def}. Henceforth we denote by $L^2(\R^2)$ and $L^2(\R_+)$ the space of square-integrable functions on $y'\in\R^2$ and $y_3\in [0,\infty)$, respectively, equipped with their usual norms. The same convention will be applied to Sobolev spaces. For the sake of notational simplicity, all norms are taken over functions on $(y',y)\in\R_+^3$ unless otherwise specified; for instance, $\|\cdot\|_{L^2}$ stands for $\|\cdot\|_{L^2(\R_+^3)}$. In addition, we also need some classical results concerning Sobolev spaces: if $\sigma>3/2$, $\tau>0$, then there are $C_1= C_1(\sigma)>0$ and $C_2= C_2(\sigma, \tau)>0$ such that
\begin{equation} \label{algebra-inq}
\|f\|_{L^\infty}\leq C_1\|f\|_{H^\sigma},\quad \|fg\|_{H^\tau}\leq C_2\|f\|_{H^\sigma}\|g\|_{H^\tau}
\end{equation}
for all $f\in H^\sigma(\R_+^3)$, $g\in H^\tau(\R_+^3)$.  Similarly, if $\sigma>1$, then  we have
\begin{equation}\label{eq:embed2d}
\|f\|_{L^\infty(\R^2)}\leq C\|f\|_{H^\sigma(\R^2)},
\end{equation}
for all $f\in H^\sigma(\R^2)$, where $C= C(\tau)>0$. In addition, if $\tau>1/2$, there is $C= C(\tau)>0$ such that the trace inequality 
\begin{equation} \label{trace}
\|f(\cdot,0)\|_{H^{\tau-1/2}(\R^2)}\leq C\|f\|_{H^\tau(\R_+^3)}
\end{equation}
holds for all $f=f(y',y_3)\in H^\tau(\R_+^3)$.

\begin{lemma}\label{le:wprop} Assume that \eqref{eq:gL} holds. Let $s, \alpha \in (5/2, \infty)$, and $\bB$ be the matrix defined in \eqref{eq:bBy}. Also, let  $\tilde{w}^\delta$ and $w$ be as in \eqref{eq:wdeltadef} and \eqref{eq:w-def}, respectively. The following statements hold.
\begin{itemize}
\item[\textup{(a)}] $\dive(\bB w)=0$.
\item[\textup{(b)}] $(\bB w)_3(y',0)=0$ for all $y'\in\R^2$, or, equivalently,
\begin{equation}\label{eq:w-lemma-eq}
w_3(y',0)=\delta^{\alpha-1}\left(D_1g(y'/\delta) w_1(y',0)+D_2g(y'/\delta) w_2(y',0)\right).
\end{equation}
\end{itemize}
In addition, there exists a positive constant $C= C(s, L)>0$ such that the following estimates hold for all $\delta \in (0,1)$.
\begin{itemize}
\item[\textup{(c)}]
\[
\|\partial_t \twd\|_{L^2}\leq C\left(\|w^0\|_{H^s}^2+\|F^0\|_{H^{s-1}}\right),
\]
and
\[
\|\partial_t w^0\|_{H^1}\leq C\left(\|w^0\|_{H^s}^2+\|F^0\|_{H^{s-1}}\right).
\]

\item[\textup{(d)}] $\|w(t,\cdot,0)\|_{H^1(\R^2)}\leq C\|w^0\|_{H^s}$.

\item[\textup{(e)}] $\|\twd\|_{H^1}\leq C\|w^0\|_{H^s}$ and $\|w\|_{H^1}\leq C\|w^0\|_{H^s}$.

\item[\textup{(f)}] $\|w\|_{L^\infty}\leq C\|w^0\|_{H^s}$ and $\|\nabla w\|_{L^\infty}\leq C\|w^0\|_{H^s}$.

\item[\textup{(g)}] $\|w(t,\cdot,0)\|_{L^\infty(\R^2)}\leq C\|w^0\|_{H^s}$.
\end{itemize}
\end{lemma}
\begin{proof}
Using \eqref{eq:bBexp}, \eqref{eq:w-def}, and that $\dive w^0=0$, a computation yields
\begin{align*}
\dive(\bB w)&=\dive\big((\bI+\delta^{\alpha-1}\bB_g)(w^0-\delta^{\alpha-5/2}\twd)\big)\\
&=\dive w^0+\delta^{\alpha-1}\dive(\bB_g w^0)-\delta^{\alpha-5/2}\dive\twd-\delta^{2\alpha-7/2}\dive(\bB_g\twd)\\
&=0+\delta^{\alpha-1}\left(-D_1g(y'/\delta) D_3w_1^0-D_2g(y'/\delta) D_3w_2^0\right)\\
&\qquad {}-\delta^{\alpha-5/2}(\delta^{3/2}+\delta^{\alpha-1})\left(-D_1g(y'/\delta) D_3w_1^0-D_2g(y'/\delta) D_3w_2^0\right)\\
&\qquad {}-\delta^{2\alpha-7/2}\left(D_1g(y'/\delta) D_3w_1^0+D_2g(y'/\delta) D_3w_2^0\right)\\
&=0,
\end{align*}
which proves statement (a).

\smallskip
In order to prove statement (b), we note that when $y_3=0$, due to the boundary condition in \eqref{eq:Euler}, we have $w_3^0(y', 0) =0$, and therefore
\[ \twd_3(y',0)=(\delta^{3/2}+\delta^{\alpha-1})\big(-D_1g(y'/\delta)w_1^0(y',0)-D_2g(y'/\delta)w_2^0(y',0)\big). 
\]
Thus,
\begin{align*}
(\bB \twd)_3(y',0)&=\big((\bI+\delta^{\alpha-1}\bB_g)\twd(y',0)\big)_3\\
&=\twd_3(y',0)+\delta^{\alpha-1}(D_1g(y'/\delta)w_1^0(y',0)+D_2g(y'/\delta)w_2^0(y',0))\\
&=\delta^{3/2} (-D_1g(y'/\delta)w_1^0(y',0)-D_2g(y'/\delta)w_2^0(y',0))\\
&=\delta^{{3/2}-\alpha+1}\delta^{\alpha-1}(-D_1g(y'/\delta)w_1^0(y',0)-D_2g(y'/\delta)w_2^0(y',0))\\
&=\delta^{5/2-\alpha}(\bB w^0)_3(y',0).
\end{align*}
The first part of statement (b) is a direct consequence of the above equality and \eqref{eq:w-def}. Equation \eqref{eq:w-lemma-eq} is obtained by expanding $(\bB w)_3$ and rearranging.

\smallskip
For statement (c), by the PDE in \eqref{eq:Euler}, we have $\partial_t w^0=\sP(-(w^0\cdot \nabla )w^0 +F^0)$, where $\sP$ is the Helmholtz-Leray projection as in Section \ref{def-weak-sol}. Applying \eqref{Projection} with $\sigma=0$ yields
\begin{align*}
\|\partial_t w^0\|_{L^2}&\leq C(\|(w^0\cdot\nabla)w^0\|_{L^2}+\|F^0\|_{L^2})\\ 
&\leq C(\|w^0\|_{H^s}^2+\|F^0\|_{H^{s-1}}),
\end{align*}
where we have used the Banach algebra property of $H^s(\R_+^3)$ stated in \eqref{algebra-inq}, and $s>5/2$. Now, $\partial_t \twd_1=-\partial_t w_1^0$ and $\partial_t \twd_2=-\partial_t w_2^0$, while
\begin{equation}\label{eq:twd3}
\partial_t \twd_3=-\partial_t w_3^0+(\delta^{3/2}+\delta^{\alpha-1})\left( -D_1g(y'/\delta) \partial_t w_1^0- D_2g(y'/\delta) \partial_t w_2^0\right).
\end{equation}
Taking $L^2$-norm on this expression, using the preceding inequality, the triangle inequality, and the fact that $\|Dg\|_{L^\infty(\R^2)}$ is bounded (as in \eqref{eq:gL}), we obtain the first inequality in statement (c), with $C= C(s, L)>0$. For the second inequality, we apply \eqref{Projection}  with $\sigma=1$ and the Banach algebra property \eqref{algebra-inq} to obtain
\begin{align*} 
\|\partial_t w^0\|_{H^1}&\leq C(\|(w^0\cdot\nabla)w^0\|_{H^1}+\|F^0\|_{H^1})\\ 
&\leq C(\|w^0\|_{H^s}^2+\|F^0\|_{H^{s-1}}),
\end{align*}
yielding the second inequality in statement (c). 

\smallskip
Statement (d) is proved using a similar argument: first, by the triangle inequality and $\delta<1$,
\[
\|w(t,\cdot,0)\|_{H^1(\R^2)}\leq \|w^0(t,\cdot,0)\|_{H^1(\R^2)}+\|\twd(t,\cdot,0)\|_{H^1(\R^2)}.
\]
For the second term, if $j=1,2$, we have $\|\twd_j(t,\cdot,0)\|_{H^1}\leq C\|w^0\|_{H^s}$ by \eqref{eq:wdeltadef} and the trace inequality \eqref{trace}. To estimate $\twd_3(t, \cdot, 0)$, we use \eqref{eq:gL}  and the last equation in \eqref{eq:wdeltadef} to obtain
\begin{multline*}
\|\twd_3(t,\cdot,0)\|_{H^1(\R^2)}\\
{}\leq C\left( \|w_3^0(t,\cdot,0)\|_{H^1(\R^2)}+\|w_1^0(t,\cdot,0)\|_{H^1(\R^2)}+\|w_2^0(t,\cdot,0)\|_{H^1(\R^2)}\right),
\end{multline*}
where we have used again that $\alpha-1>3/2$ to infer that $C$ is independent of $\delta \in (0,1)$. Using the trace inequality \eqref{trace} again, statement (d) is obtained. 

\smallskip
The first inequality in statement (e) is a consequence of \eqref{eq:wdeltadef}, $\alpha>5/2$, the triangle inequality, and \eqref{eq:gL}, while the second inequality is obtained directly from the first inequality, the triangle inequality, and \eqref{eq:w-def}. Statement (f) follows from \eqref{eq:w-def}, \eqref{eq:wdeltadef}, the triangle inequality, the assumption \eqref{eq:gL},  the Sobolev embedding \eqref{algebra-inq}, and the fact that $\alpha>5/2$.

\smallskip
For statement (g), from \eqref{eq:w-def} and \eqref{eq:wdeltadef}, we see that
\[
\|w(t,\cdot,0)\|_{L^\infty(\R^2)}\leq C\|w^0(t,\cdot,0)\|_{L^\infty(\R^2)}.
\]
Applying \eqref{eq:embed2d} and \eqref{trace}, statement (g) follows, and the proof of the lemma is complete.
\end{proof}
\subsection{Construction of $\sB$ and its estimates} We build the term $\sB$, the second term in the boundary layer $\hat{\sB}$ defined in \eqref{hat-sB-term}, and provide its estimates. Let us denote 
\[ v = u - w^0 - \hat{\sB}  = u - w - \sB.
\] 
Due to \eqref{diff-bdr-cond} and the discussion right after it, we require that $\sB$ has small $L^2$ norm, is localized near $\{y_3 =0\}$, and
\begin{equation} \label{B.cond}
\dive(\bB\sB) =0, \quad \sB\big|_{\{y_3 =0\}} = -w\big|_{\{y_3 =0\}}, \quad \text{and} \quad  \sB\big|_{\{y_3 = \infty\}} =0.
\end{equation}
We note that it follows from the definition of $\bB$ in \eqref{eq:bBy} that $\dive(\bB \sB) =0$ is equivalent to
\begin{equation} \label{Div.cond}
\dive{\sB} =  \delta^{\alpha-1}D_1 g(y'/\delta) D_3 \sB_1  + \delta^{\alpha-1}D_2 g(y'/\delta) D_3 \sB_2 .
\end{equation}
Following \cite{PV}, which in turn is inspired by \cite{Mas}, we define $\sB:(0,T^*)\times\R_+^3\to\R^3$ as
\begin{equation} \label{BL.def}
\begin{split}
\sB_1(t, y',y_3)&= -w_1(t, y',0)\phi(y_3/\sqrt{\theta\nu}),\\
\sB_2(t, y',y_3)&= -w_2(t, y',0)\phi(y_3/\sqrt{\theta\nu}),\\
\sB_3(t, y',y_3)&=-w_3(t, y', 0)\phi(y_3/\sqrt{\theta\nu}) \\
 & \quad  +\sqrt{\theta\nu}[D_1w_{1}(t, y',0)+ D_2w_{2}(t, y',0)]\psi(y_3/\sqrt{\theta\nu}),
\end{split}
\end{equation}
where $\theta>0$ is a small constant to be determined, $y' \in \R^2$, $y_3 >0$, $t \in (0, T^*)$, and $\phi$, $\psi$ are smooth functions on $[0,\infty)$ satisfying 
\begin{equation} \label{eq:PsiPhi}
\begin{split}
& \psi'=\phi, \quad  \phi(0)=1,  \quad  \psi(0)=0,  \quad \text{and} \\
& \phi(z)=\psi(z)\equiv 0 \quad \text{for }  z>1. 
\end{split}
\end{equation}
We note that in \eqref{BL.def} the function $w$ defined in \eqref{eq:w-def} is the perturbation of the solution $w^0$ to \eqref{eq:Euler}, and $w$ depends on $\delta$, but we omit this dependence from the notation for the sake of simplicity. The functions $\phi$ and $\psi$, on the other hand, are henceforth fixed. Note also that, unlike \cite{PV}, we use $w$ instead of $w^0$ in the construction of $\sB$.

\smallskip
It is clear from our choice of $\phi$, $\psi$ that the boundary conditions
\begin{equation*}
\sB(t, y', 0) = -w(t, y',0), \quad \sB(t, y', \infty) =0 \quad (t, y') \in (0, T^*) \times \R^2
\end{equation*}
hold. Similarly,
\[
\begin{split}
\dive \sB & = -[D_1 w_1(t, y', 0) + D_2 w_2(t, y',0)] \phi(y_3/\sqrt{\theta\nu}) \\
&\quad {}-\frac{1}{\sqrt{\theta\nu}} w_3(t, y', 0) \phi'(y_3/\sqrt{\theta\nu }) \\
& \quad {}+  [D_1w_{1}(t, y',0)+ D_2w_{2}(t, y',0)]\psi'(y_3/\sqrt{\theta\nu})\\
& =  -\frac{1}{\sqrt{\theta\nu}} w_3(t, y', 0) \phi'(y_3/\sqrt{\theta\nu}),
\end{split}
\]
and
\[
\begin{split}
& \delta^{\alpha-1} D_1 g(y'/\delta) D_3 \sB_1  + \delta^{\alpha-1} D_2 g(y'/\delta) D_3 \sB_2  \\
& = -\frac{\delta^{\alpha-1}}{\sqrt{\theta\nu}} \big[  D_1 g(y'/\delta) w_1(t, y', 0)  +   D_2 g(y'/\delta) w_2(t, y', 0) \big] \phi'(y_3/\sqrt{\theta\nu}) \\
& =-\frac{1}{\sqrt{\theta\nu}} w_3(t, y', 0) \phi'(y_3/\sqrt{\theta\nu}),
\end{split}
\]
where we used \eqref{eq:w-lemma-eq} in the last identity. Hence \eqref{Div.cond} holds, and therefore $\dive(\bB\sB)=0$.
\begin{rmkn}\label{rmk:sB}
Note that $\sB_1$, $\sB_2$, and each term in $\sB_3$ is the product of a function of $y' \in \mathbb{R}^2$ and a function of $y_3 \in (0, \infty)$. This is useful in some calculations. For example, we have
\[ \|\sB_1(t,\cdot)\|_{L^2}=\|w_1(t,\cdot,0)\|_{L^2(\R^2)}\|\phi(\cdot/\sqrt{\theta\nu})\|_{L^2(\R_+)}.
\]
\end{rmkn}

\smallskip
Next, in Proposition \ref{prop:sBbounds} below, we derive essential estimates  of $\sB$ that are needed in the next section. For simplicity, let us denote
\begin{equation}\label{eq:sB}
\sB(t, y) =- w(t, y', 0) \phi(y_3/\sqrt{\theta\nu}) + \sqrt{\theta\nu} \vec{e}_3\psi(y_3/\sqrt{\theta\nu}) A(t, y'), 
\end{equation}
for $y = (y', y_3) \in \R^2 \times \R_+$ and $t \in (0, T^*)$,
where
\begin{equation}\label{eq:Aty}
A(t, y') = D_1w_{1}(t, y',0)+ D_2w_{2}(t, y',0), \quad \vec{e}_3 = (0,0,1).
\end{equation}
Since $A$ involves only the first and second components of $w$, which do not involve $Dg$ (see \eqref{eq:wdeltadef} and \eqref{eq:w-def}), one can use the triangle inequality, the trace inequality \eqref{trace}, and the Sobolev embedding \eqref{eq:embed2d} to obtain
\begin{equation}\label{eq:Abounds}
\|A(t,\cdot)\|_{H^1(\R^2)}\leq C\|w^0\|_{H^s}\quad\text{and}\quad \|A(t,\cdot)\|_{L^\infty(\R^2)}\leq C\|w^0\|_{H^s},
\end{equation}
where $C= C(s)>0$ is a constant independent of $\delta \in (0,1)$. These inequalities will be useful later on.

\smallskip
For the rest of this section, although the functions $\phi$, $\psi$ are functions of a single variable, for convenience we write $D_3[\phi]$, $D_3[\psi]$ to denote the derivative operator acting on them as they are mainly used with the $y_3$-variable. Moreover, for each $p \in [1, \infty)$ and $\gamma \in \R$, we denote the weighted Lebesgue space $L^p(\R_+, y_3^\gamma)$ with norm
\[
\|f\|_{L^p(\R_+, y_3^\gamma)} = \left(\int_0^\infty |f(y_3)|^p y_3^\gamma dy_3 \right)^{1/p},
\]
and a similar definition can be formulated when $p=\infty$. 

\smallskip
To prove Proposition \ref{prop:sBbounds} below, we need the following lemma giving some simple $L^2$ and $L^\infty$ weighted and unweighted estimates of $\phi$ and $\psi$.
\begin{lemma}\label{le:bounds1} There exists a constant $C$ depending on $\phi$ and $\psi$, such that the following estimates hold for all $\theta, \nu \in (0,1)$:
\begin{itemize}
\item[\textup{(a)}] $\|\phi(\cdot/\sqrt{\theta\nu})\|_{L^2(\R_+)}\leq C(\theta\nu)^{1/4}$,
\item[\textup{(b)}] $\|\sqrt{\theta\nu}\psi(\cdot/\sqrt{\theta\nu})\|_{L^2(\R_+)}\leq C(\theta\nu)^{1/4}$,
\item[\textup{(c)}] $\|D_3[\phi(\cdot/\sqrt{\theta\nu})]\|_{L^2(\R_+, y_3^2)}\leq C(\theta\nu)^{1/4}$,
\item[\textup{(d)}] $\|D_3[\phi(\cdot /\sqrt{\theta\nu})]\|_{L^2(\R_+)}\leq C(\theta\nu)^{-1/4}$,
\item[\textup{(e)}] $\|D_3[\psi(\cdot /\sqrt{\theta\nu})]\|_{L^2(\R_+, y_3^2)}\leq C(\theta\nu)^{1/4}$,
\item[\textup{(f)}] $\|\sqrt{\theta\nu}\psi(\cdot/\sqrt{\theta\nu})\|_{L^\infty(\R_+)}\leq C\sqrt{\theta\nu}$,
\item[\textup{(g)}] $\|D_3[\phi(\cdot /\sqrt{\theta\nu})]\|_{L^\infty(\R_+, y_3^2)}\leq C\sqrt{\theta\nu}$,
\item[\textup{(h)}] $\|D_3[\psi(\cdot/\sqrt{\theta\nu})]\|_{L^\infty(\R_+, y_3^2)}\leq C\sqrt{\theta\nu}$.
\end{itemize}
\end{lemma}
\noindent
The proof of this lemma is elementary, and we skip it. Readers can find the details of its proof in \cite[Lemma 2, Section 3]{PV}.

\smallskip
In the statement of the Proposition \ref{prop:sBbounds} below, to further simplify the notation, we omit the arguments $(t,\cdot)$ from the functions $\sB$, $w^0$, and $F^0$. In particular, all the forthcoming estimates are valid for all $t\in (0,T^*)$, where $T^*$ is the time such that the strong solution $w^0$ of \eqref{eq:Euler} exists on $(0, T^*)$. This proposition is analogous to \cite[Proposition 1, Section 3]{PV}.
\begin{prop}\label{prop:sBbounds} Let $s$, $\alpha \in (5/2, \infty)$ and assume that \eqref{eq:gL}  holds. There exists a constant $C>0$ depending on $\phi$, $\psi, s$, and $L$ such that the following estimates hold for each $t \in (0, T^*)$, $\theta, \delta, \nu \in (0,1)$\textup{:}
\begin{itemize}
\item[\textup{(a)}] $\|\partial_t\sB\|_{L^2(\R^3_+)}\leq C(\theta\nu)^{1/4}\big(\|w^0\|_{H^s(\R^3_+)}^2+\|F^0\|_{H^{s-1}(\R^3_+)}\big)$,
\item[\textup{(b)}] $\|D_j\sB\|_{L^2(\R^3_+)}\leq C(\theta\nu)^{1/4}\|w^0\|_{H^s(\R^3_+)}$ \text{for} $j \in \{1,2\}$,
\item[\textup{(c)}] $\|D_3\sB\|_{L^2(\R^3_+, y_3^2)}\leq C(\theta\nu)^{1/4}\|w^0\|_{H^s(\R^3_+)}$ and\\ $\|D_3\sB\|_{L^2(\R^3_+)}\leq C(\theta\nu)^{-1/4}\|w^0\|_{H^s}$,
\item[\textup{(d)}] $\|D_3\sB\|_{L^\infty(\R^3_+, y_3^2)}\leq C\sqrt{\theta\nu}\|w^0\|_{H^s(\R^3_+)}$,
\item[\textup{(e)}] $\|\bB\sB\|_{L^2(\R^3_+)}\leq C(\theta\nu)^{1/4}\|w^0\|_{H^s(\R^3_+)}$ and\\ $\|[\bB\sB]_3\|_{L^2(\R^3_+)}\leq C\sqrt{\theta\nu}\|w^0\|_{H^s(\R^3_+)}$,
\item[\textup{(f)}] $\|[\bB\sB]_j\|_{L^\infty(\R^3_+)}\leq C\|w^0\|_{H^s(\R^3_+)}$ \text{for} $j \in \{1,2\}$, and\\ $\|[\bB\sB]_3\|_{L^\infty(\R^3_+)}\leq C\sqrt{\theta\nu}\|w^0\|_{H^s(\R^3_+)}$.
\end{itemize}
\end{prop}
\begin{proof} 
The proofs of statements (c), (d), and (f) are analogous to the corresponding statements in \cite[Proposition 1, Section 3]{PV}, using Lemma \ref{le:wprop} (d), (f), and (g) when necessary to obtain bounds depending on $w^0$ rather than on $w$, as well as \eqref{eq:Abounds} for bounds involving $A$ as in \eqref{eq:Aty}. Therefore, we only provide detailed proofs of statements (a), (b), and (e).

\smallskip
\noindent To prove statement (a), we first note that, by the definitions of $\twd$ and $w$ in \eqref{eq:wdeltadef}--\eqref{eq:w-def}, distributing the derivative with respect to $t$, one has
\begin{align*}
\|\partial_t w(t,\cdot,0)\|_{L^2(\R^2)}&\leq \|\partial_t w^0(t,\cdot,0)\|_{L^2(\R^2)}+\|\partial_t \twd(t,\cdot,0)\|_{L^2(\R^2)}\\
&\leq C\|\partial_t w^0(t,\cdot,0)\|_{L^2(\R^2)} \\
&\leq C\|\partial_t w^0\|_{H^1}\\
&\leq C\left(\|w^0\|_{H^s}^2+\|F^0\|_{H^{s-1}}\right),
\end{align*}
using the trace inequality \eqref{trace} (with $\tau=1/2$) and Lemma \ref{le:wprop} (c). A similar procedure applied to $A$ in \eqref{eq:Aty} implies that
\[
\|\partial_t A(t,\cdot)\|_{L^2(\R^2)}\leq C\left(\|w^0\|_{H^s}^2+\|F^0\|_{H^{s-1}}\right).
\]
Using \eqref{eq:sB}, \eqref{eq:Aty}, the product form of $\sB$ as in Remark \ref{rmk:sB}, and Lemma \ref{le:bounds1} (a) and (b), we obtain
\begin{align*}
\|\partial_t\sB\|_{L^2}&\leq \|\partial_t w(t,\cdot,0)\|_{L^2(\R^2)} \|\phi(\cdot/\sqrt{\theta\nu})\|_{L^2(\R_+)}\\
&\qquad{}+\sqrt{\theta\nu}\|\psi(\cdot/\sqrt{\theta\nu})\|_{L^2(\R_+)}\|\partial_t A(t,\cdot)\|_{L^2(\R^2)}\\
&\leq C\bigg(\|\partial_t w(t,\cdot,0)\|_{L^2(\R^2)} (\theta\nu)^{1/4}+(\theta\nu)^{1/4}\|\partial_t A(t,\cdot)\|_{L^2(\R^2)}\bigg)\\
&\leq C(\theta\nu)^{1/4}\left(\|w^0\|_{H^s}^2+\|F^0\|_{H^{s-1}}\right),
\end{align*}
where we have used the foregoing estimates on $\|\partial_t w(t,\cdot,0)\|_{L^2(\R^2)}$ and $\|\partial_t A(t,\cdot)\|_{L^2(\R^2)}$.

\smallskip
\noindent Next, we prove (b). From Lemma \ref{le:bounds1} (a) and (b), Remark \ref{rmk:sB}, and \eqref{eq:Abounds},
\begin{align*}
\|D_j\sB\|_{L^2}&\leq \|D_jw(t,\cdot,0)\|_{L^2(\R^2)}\|\phi(\cdot/\sqrt{\theta\nu})\|_{L^2(\R_+)}\\
&\qquad{}+\|\sqrt{\theta\nu}\psi(\cdot/\sqrt{\theta\nu})\|_{L^2(\R_+)} \|D_jA(t,\cdot)\|_{L^2(\R^2)}\\
&\leq C(\theta\nu)^{1/4}\left(\|w(t,\cdot,0)\|_{H^1(\R^2)}+\|w^0\|_{H^s}\right)\\
&\leq C(\theta\nu)^{1/4}\|w^0\|_{H^s},
\end{align*}
the last inequality making use of Lemma \ref{le:wprop} (d), while the second to last inequality uses \eqref{eq:Abounds}.

\smallskip
\noindent Now, we prove (e). We first note that for $j\in\{1,2\}$ we have $[\bB\sB]_j=\sB_j$, so \eqref{eq:sB}, the product form of $\sB$ as in Remark \ref{rmk:sB}, and Lemma \ref{le:bounds1} (a) imply 
\begin{align*}
\|[\bB\sB]_j\|_{L^2}&=\|w(t,\cdot,0)\|_{L^2(\R^2)}\|\phi(\cdot/\sqrt{\theta\nu})\|_{L^2(\R_+)}\\
&\leq C\|w^0\|_{H^s}(\theta\nu)^{1/4}.
\end{align*}
On the other hand, we can rewrite $[\bB\sB]_3$ as follows:
\begin{align*}
[\bB\sB]_3&=-\delta^{\alpha-1}D_1g(y')\sB_1(t,y)-\delta^{\alpha-1}D_2g(y')\sB_2(t,y)+\sB_3(t,y)\\
&=\delta^{\alpha-1}D_1g(y')w_1(t,y',0)\phi(y_3/\sqrt{\theta\nu})\\
&\qquad {}+\delta^{\alpha-1}D_2g(y')w_2(t,y',0)\phi(y_3/\sqrt{\theta\nu})-w_3(t,y',0)\phi(y_3/\sqrt{\theta\nu})\\
&\qquad {}+\sqrt{\theta\nu}[D_1w_1(t,y',0)+D_2w_2(t,y',0)]\psi(y_3/\sqrt{\theta\nu})\\
&=\sqrt{\theta\nu}A(t,y')\psi(y_3/\sqrt{\theta\nu}),
\end{align*}
where we have used that the first three terms are cancelled by \eqref{eq:w-lemma-eq}. From \eqref{eq:Abounds}, Remark \ref{rmk:sB}, Lemma \ref{le:bounds1} (f), and Lemma \ref{le:wprop} (d),
\begin{align*}
\|[\bB\sB]_3\|_{L^2}&\leq \|\sqrt{\theta\nu}\psi(\cdot/\sqrt{\theta\nu})\|_{L^\infty(\R_+)}\|A(t, \cdot)\|_{L^2(\R^2)}\\
&\leq C\sqrt{\theta\nu}\|w^0\|_{H^s}.
\end{align*}
The proof of the proposition is thus complete.
\end{proof}

\section{Energy estimates}\label{Ener-sec}  Let $T^* \in (0, T)$ be fixed such that the strong solution 
\[ w^0 \in L^\infty((0, T^*), H^s(\R_+^3)^3)\]
of \eqref{eq:Euler} exists as in Lemma \ref{Euler-existence}. For the given $w^0$, let $w$ be defined in \eqref{eq:w-def} and $\sB$ be defined in \eqref{BL.def}. Also, let  $u(t,y) = \bar{u}(t, \Psi_0(y))$, where $\bar{u}$ is a Leray-Hopf weak solution of \eqref{eq:NSE} defined as in Subsection \ref{def-weak-sol}. Observe that $u$ is defined on $(0, T^*) \times \R^3_+$ due to the change of coordinates as in \eqref{change-variables}, and $u$ is a weak solution to \eqref{NS.eqn}. 

\smallskip
Let us denote
\[ 
v = u -w - \sB. 
\]
The following proposition on the energy estimate of $v$ is the main result of the section.
\begin{prop}\label{prop:G} Let $\Lambda \in (0,1)$, $L$, $K_0$ be positive numbers, $\alpha$, $s \in (5/2, \infty)$, and let $\delta_0 = \delta_0 (\Lambda, L, K_0)>0$ be defined in \textup{Proposition \ref{prop:Z}}. Assume that the assumptions in \textup{Theorem \ref{thm:main}} hold. Then, there exist nonnegative functions $f_0(t)$, $f_1(t)$, and $f_2(t)\in L^1(0,T^*)$ such that the following inequality holds:
\begin{equation}\label{eq:ineq}
\begin{split}
& \frac{1}{2}\partial_t\|v\|_{L^2}^2+\int_{\R^3_+} \wei{\A_0 \bB^* \nabla v, \bB^* \nabla v}\, dy\\
& \quad\leq f_0(t)\|v\|_{L^2}^2+ f_1(t)\|v\|_{L^2}+f_2(t), \quad t \in (0, T^*).
\end{split}
\end{equation}
Moreover, the functions $f_0$, $f_1$, and $f_2$ satisfy
\begin{equation}\label{eq:abcboundsbeta}
\begin{split}
& \int_0^{T^*} f_0(t)\,dt\leq C ,\quad \int_0^{T^*} f_1(t)\,dt\leq C(\beta+\delta^{\alpha-5/2}) ,\\
& \int_0^{T^*} f_2(t)\,dt\leq C(\beta+\delta^{\alpha-5/2})^2,
\end{split}
\end{equation}
where $C$ is a constant depending on $\|w^0\|_{L^\infty((0,T^*),  H^s)}$, $\|F^0\|_{L^1((0,T^*),H^s(\R_+^3))}$, $s$, $T^*$, $L$, $K_0$, and $\Lambda$, but independent of $u$, $\eta$, $\nu$, $\delta$, and $\beta$, where $\beta$ is the function defined in \textup{Definition \ref{hyp:thm}}.
\end{prop}
\smallskip
The proof of Proposition \ref{prop:G} is similar to that of \cite[Proposition 2, Section 3]{PV}, with changes to account for the setting of the present paper. The changes needed involve both the new terms arising from the specifics of the problem we consider in this article, as well as from the fact that the ellipticity assumptions used in \cite{PV} (see Hypothesis 1 (ii) in \cite{PV}) do not directly apply in our context; instead, estimates involving ellipticity will make use of Proposition \ref{prop:Z}. 

\smallskip
To proceed, we derive the equation for $v$. Note that $\dive(\bB v)=0$ by \eqref{NS.eqn}, \eqref{B.cond}, and Lemma \ref{le:wprop} (a). In addition, as $u = v + w + \sB$, it follows from \eqref{NS.eqn}  that
\begin{equation}\label{eq:pretest}
\begin{aligned}
\partial_t v+ \partial_tw & + \partial_t \sB  - \dive[\A \nabla ( v+ w+ \sB)]\\
& + [(\bB ( v +w +\sB))\cdot \nabla] (v+ w+ \sB)
 = F - \bB^* \nabla p,
\end{aligned}
\end{equation}
where $\A$ is the matrix in \eqref{eq:bA}. Due to the form of $w$ in \eqref{eq:w-def} and the decomposition of $\bB$ in \eqref{eq:bBexp}, we see that
\begin{align*}
[(\bB w)\cdot\nabla]w &=\bigg[ \big((\bI+\delta^{\alpha-1}\bB_g)w\big)\cdot\nabla\bigg]w\\
&=w\cdot\nabla w+\delta^{\alpha-1}\big[(\bB_gw)\cdot\nabla\big]w\\
&=\big[ \big(w^0-\delta^{\alpha-5/2} \twd\big)\cdot\nabla\big](w^0-\delta^{\alpha-5/2}\twd)\\
&\qquad {}+\delta^{\alpha-1}\big[(\bB_gw)\cdot\nabla\big]w\\
&=(w^0\cdot\nabla)w^0-\delta^{\alpha-5/2}(w^0\cdot\nabla)\twd-\delta^{\alpha-5/2}(\twd\cdot\nabla)w\\
&\qquad {}+\delta^{\alpha-1}\big[(\bB_gw)\cdot\nabla\big]w.
\end{align*}
Using this equality, \eqref{eq:w-def}, and the fact that $w^0$ satisfies \eqref{eq:Euler}, we rewrite equation \eqref{eq:pretest} as
\begin{align*}
\partial_t v & -\delta^{\alpha-5/2}\partial_t \twd  +\partial_t\sB-\dive[\A\nabla(v+w+\sB)]\\
&\qquad  +[\bB w\cdot\nabla](v+\sB)+[\bB(v+\sB)\cdot\nabla](v+w+\sB)\\
& \qquad -\delta^{\alpha-5/2}(w^0\cdot\nabla) \twd-\delta^{\alpha-5/2}(\twd\cdot\nabla)w + \delta^{\alpha-1}[(\bB_g w)\cdot\nabla]w\\
& =F-F^0-\bB^*\nabla p+\nabla q \qquad \text{in }(0, T^*) \times \mathbb{R}^3_+.
\end{align*}
\smallskip
Now, testing this PDE with $v$, using the divergence free condition $\dive(\bB v) =0$, the boundary condition $v\big|_{\{y_3=0\}}$, and rearranging, we obtain
\[
\begin{split}
& \frac{1}{2}\frac{d}{dt} \int_{\R^3_+} |v|^2\, dy + \int_{\R^3_+} \wei{\A \nabla v, \nabla v}\, dy = \int_{\R^3_+} \wei{- \partial_t\sB + F-F^0,v}\, dy \\
& {}-\int_{\R^3_+} \big(\wei{\A ( \nabla w+ \nabla \sB), \nabla v} + \wei{([\bB w]\cdot \nabla )(v+ \sB), v} \big)\, dy \\
& {}- \int_{\R^3_+}  \wei{ ([\bB (v+\sB)]\cdot \nabla )(v+ w+\sB), v}\, dy\\
& {}+\delta^{\alpha-5/2}\bigg(\int_{\R_+^3}\wei{\partial_t \twd,v}\,dy+\int_{\R_+^3}\wei{(w^0\cdot\nabla)\twd,v}\,dy+\int_{\R_+^3}\wei{(\twd\cdot\nabla)w,v}\,dy\bigg)\\
& {}-\delta^{\alpha-1}\int_{\R_+^3}\wei{[(\bB_g w)\cdot\nabla]w,v}\,dy+\int_{\R_+^3}\wei{\nabla q,v}\,dy.
\end{split}
\]
Here, for notational clarification, we note that
\[
\wei{\A \nabla v, \nabla v} =\sum_{k=1}^3 \wei{\A\nabla v_k,\nabla v_k}, \quad  \wei{[\bB w]\cdot \nabla v, v}  = \sum_{i,j,k=1}^3  v_k\bB_{ij} w_jD_iv_k,
\]
and similarly for other terms. Rearranging terms, and recalling \eqref{eq:bA} that $\A=\bB\A_0\bB^*$ where $\bB$ is the matrix in \eqref{eq:bBy}, we then obtain  
\begin{align} \notag
& \frac{1}{2}\frac{d}{dt}\int_{\R^3_+} |v|^2\, dy + \int_{\R^3_+} \wei{\A_0 \bB^* \nabla v, \bB^* \nabla v}\, dy \\ \notag
= &  -\int_{\R^3_+} \big[ \wei{([\bB w]\cdot\nabla)\sB,v}  + \wei{([\bB\sB]\cdot\nabla)\sB,v}  + \wei{([\bB\sB]\cdot\nabla)w,v}\big]\,dy\\ \notag
& -\int_{\R^3_+} \wei{([\bB v]\cdot\nabla)w,v}\,dy -\int_{\R^3_+} \big[ \wei{\partial_t\sB,v}+\wei{([\bB(v+w+\sB)]\cdot\nabla)v,v}\big]\,dy\\ \notag
& -\int_{\R^3_+} \wei{\A ( \nabla w+ \nabla \sB), \nabla v}\,dy-\int_{\R^3_+} \wei{([\bB v]\cdot\nabla)\sB,v}\,dy+\int_{\R^3_+} \wei{F-F^0,v}\,dy\\
& {}+\delta^{\alpha-5/2}\bigg(\int_{\R_+^3}\wei{\partial_t \twd,v}\,dy+\int_{\R_+^3}\wei{(w^0\cdot\nabla)\twd,v}\,dy+\int_{\R_+^3}\wei{(\twd\cdot\nabla)w,v}\,dy\bigg)\nonumber\\
& {}-\delta^{\alpha-1}\int_{\R_+^3}\wei{[(\bB_g w)\cdot\nabla]w,v}\,dy+\int_{\R_+^3}\wei{\nabla q,v}\,dy. \label{eq:Mas17}
\end{align}

\smallskip
We now provide a series of lemmas in which we obtain certain bounds on each integral on the right hand side of \eqref{eq:Mas17}. We also recall that in Theorem \ref{thm:main}, we assume 
\[
s, \alpha \in (5/2, \infty),\  \eqref{eq:gL}, \ \eqref{ellipticity-cond},  \text{ and } \delta^{\alpha-1} \eta \leq K_0 \nu.
\]
Unless otherwise specified, in the forthcoming lemmas the constant $C$ may depend on $L$, $\Lambda$, $K_0$, and $s$, but is independent of $u$, $w$, $w^0$, $\eta$, $\nu$, and $\delta$. We first state three lemmas providing estimates of the first, second, and third lines on the right hand side of \eqref{eq:Mas17}.
\begin{lemma} \label{123-est.0816} Let  $s$, $\alpha \in (5/2, \infty)$, and assume that \eqref{eq:gL} holds. There is a constant $C= C(L,s)>0$ such that the following estimates
\begin{align*}
\textup{(i)} \quad  & \left|\int_{\R^3_+} \wei{([\bB w]\cdot\nabla)\sB,v}\,dy\right|+\left|\int_{\R^3_+} \wei{([\bB\sB]\cdot\nabla)\sB,v}\,dy\right| \\
& \quad {}+\left|\int_{\R^3_+} \wei{([\bB\sB]\cdot\nabla)w,v}\,dy\right|\leq C(\theta\nu)^{1/4}\|w^0\|_{H^s}^2 \|v\|_{L^2},\\
\textup{(ii)} \quad  & \left|\int_{\R^3_+} \wei{([\bB v]\cdot\nabla)w,v}\,dy\right|\leq C\|w^0\|_{H^s}\|v\|_{L^2}^2, \\
\textup{(iii)} \quad  & \left|\int_{\R^3_+} \wei{\partial_t\sB,v}\,dy\right|+\left|\int_{\R^3_+}\wei{([\bB(v+w+\sB)]\cdot\nabla)v,v}\,dy\right|\\ 
& \qquad \leq C(\theta\nu)^{1/4}\left(\|w^0\|_{H^s}^2+\|F^0\|_{L^{2}}\right)\|v\|_{L^2}
\end{align*}
hold for all positive constants $\theta, \nu, \delta \in (0,1)$.
\end{lemma}
\begin{proof} The proof of (i) and (ii) follows by the direct application of H\"{o}lder's inequality using Proposition \ref{prop:sBbounds} to control terms involving $\sB$, and Lemma \ref{le:wprop} (f) to control terms involving $w$. As the calculations are straightforward and they are similar to the proofs of \cite[Lemmas 3, 4]{PV}, we skip the details.
\smallskip

To prove (iii), we note that $v+w +\sB =u$ and 
\[
\int_{\R^3_+}\wei{([\bB(v+w+\sB)]\cdot\nabla)v,v}\,dy =\int_{\R^3_+}\wei{(\bB u\cdot\nabla)v,v}\,dy =0 
\]
by using integration by parts, and the fact that $\text{div}(\bB u)=0$. On the other hand, by Proposition \ref{prop:sBbounds} (a), we have
\[
\left|\int_{\R^3_+} \wei{\partial_t\sB,v}\,dy\right| \leq \|\partial_t\sB\|_{L^2} \|v\|_{L^2} \leq C(\theta\nu)^{1/4}\left(\|w^0\|_{H^s}^2+\|F^0\|_{L^{2}}\right)\|v\|_{L^2}.
\]
Therefore (iii) follows and the proof is then completed.
\end{proof}

\smallskip
We recall the notation $D_{y'}=(D_1,D_2)$. The fourth integral on the right hand side of \eqref{eq:Mas17} is now controlled in the following lemma.

\begin{lemma}\label{le:Mas811}  Let  $s, \alpha \in (5/2, \infty)$, $K_0>0$, and assume that \eqref{eq:gL} and  \eqref{ellipticity-cond} hold. There exists a constant $C= C(\Lambda, L, K_0, s)>0$ such that
\begin{align*}
& \left|\int_{\R^3_+} \wei{\A ( \nabla w+ \nabla \sB), \nabla v}\,dy\right|\leq  \varepsilon \int_{\R^3_+} \wei{\A_0 \bB^*\nabla v,\bB^*\nabla v}\,dy\\
& \qquad {}+\frac{C}{\varepsilon}\left[\eta(1+\sqrt{\theta\nu})+\nu(1+(\theta\nu)^{-1/2})\right]\|w^0\|_{H^s}^2
\end{align*}
holds for any $\varepsilon>0$, and $\theta, \eta \in (0,1), \nu \in (0, \eta)$, and for $\delta \in (0, \delta_0)$ satisfying $\delta^{\alpha-1} \eta \leq K_0 \nu$.
\end{lemma}
\begin{proof} From \eqref{eq:bA} and Proposition \ref{prop:Z} (b), we have  
\begin{align*}
& |\wei{\A ( \nabla w+ \nabla \sB), \nabla v}|=|\wei{\A_0\bB^*( \nabla w+ \nabla \sB),  \bB^*\nabla v}|\\
& \leq C\Big[ \eta|D_{y'} w+ D_{y'} \sB || D_{y'} v|+ \nu|D_3 w+ D_3 \sB| |D_3 v|\\
&\qquad {}+(\delta^{\alpha-1}\eta+\nu)\big(|D_3w+D_3\sB ||D_{y'}v|+|D_3v| |D_{y'}w+D_{y'}\sB|\big) \Big],
\end{align*}
where $C = C(\Lambda, L, K_0) >0$. Then
\begin{equation}\label{eq:Awbv}
\begin{aligned}
& \left|\int_{\R^3_+} \wei{\A ( \nabla w+ \nabla \sB), \nabla v}\,dy\right| \\
&\leq C \bigg[ \eta \int_{\R^3_+} |D_{y'} w+ D_{y'}\sB || D_{y'} v|\,dy + \nu \int_{\R^3_+} |D_3 w+ D_3\sB || D_3 v|\,dy\\
&\qquad {}+(\delta^{\alpha-1}\eta+\nu)\int_{\R_+^3} |D_3w+D_3\sB||D_{y'}v|\,dy\\
&\qquad {}+(\delta^{\alpha-1}\eta+\nu)\int_{\R_+^3}|D_3v||D_{y'}w+D_{y'}\sB|\,dy \bigg].
\end{aligned}
\end{equation}
Recalling from Proposition \ref{prop:sBbounds} (b) and (c) that 
\begin{equation} \label{add-Nov-24-24}
\|D_{y'}\sB\|_{L^2} \leq C(\theta\nu)^{1/4}\|w^0\|_{H^s},\quad  \|D_3\sB\|_{L^2} \leq C(\theta\nu)^{-1/4}\|w^0\|_{H^s},
\end{equation}
and noting that $\|Dw\|_{L^2}\leq C\|w^0\|_{H^s}$ by Lemma \ref{le:wprop} (e), we find
\begin{multline*}
\eta \int_{\R^3_+} |D_{y'} w+ D_{y'}\sB || D_{y'} v|\,dy + \nu \int_{\R^3_+} |D_3 w+ D_3\sB || D_3 v|\,dy \\
\leq C\Big[\eta \big( 1+(\theta\nu)^{1/4}\big)\|w^0\|_{H^s}\|D_{y'} v\|_{L^2}+ \nu \big(1+(\theta\nu)^{-1/4}\big)\|w^0\|_{H^s}\|D_3 v\|_{L^2}\Big].
\end{multline*}
For the other two integrals, from H\"{o}lder's inequality and \eqref{add-Nov-24-24}, we have
\begin{align*}
\int_{\R_+^3} |D_3w+D_3\sB||D_{y'}v|\,dy&\leq C\|D_3w+D_3\sB\|_{L^2}\|D_{y'}v\|_{L^2}\\
&\leq C\big( \|w^0\|_{H^s}+(\theta\nu)^{-1/4}\|w^0\|_{H^s}\big)\|D_{y'}v\|_{L^2}\\
&=C\|D_{y'}v\|_{L^2} \|w^0\|_{H^s}\big(1+(\theta\nu)^{-1/4}\big),
\end{align*}
and
\begin{align*}
\int_{\R_+^3} |D_3v||D_{y'}w+D_{y'}\sB|\,dy&\leq C\|D_3v\|_{L^2}\|D_{y'}w+D_{y'}\sB\|_{L^2}\\
&\leq C\|D_3v\|_{L^2}\big(\|w^0\|_{H^s}+(\theta\nu)^{1/4}\|w^0\|_{H^s}\big)\\
&= C\|w^0\|_{H^s}\|D_3 v\|_{L^2}\big(1+(\theta\nu)^{1/4}\big).
\end{align*}
Substitute the last three inequalities into \eqref{eq:Awbv} and rearrange to find
\begin{align*}
&\left|\int_{\R_+^3}\wei{\A(\nabla w+\nabla\sB,\nabla v}\,dy\right|\\
&\leq C\Big[ \big(\eta+\eta(\theta\nu)^{1/4}+(\delta^{\alpha-1}\eta+\nu)(1+(\theta\nu)^{-1/4})\big)\|w^0\|_{H^s}\|D_{y'} v\|_{L^2}\\
&\qquad{}+ \big(\nu+\nu(\theta\nu)^{-1/4}+(\delta^{\alpha-1}\eta+\nu)(1+(\theta\nu)^{1/4})\big)\|w^0\|_{H^s}\|D_3 v\|_{L^2}\Big]\\
&\leq C\Big[ \big( \eta+\eta(\theta\nu)^{1/4}+\nu(\theta\nu)^{-1/4}\big)\|w^0\|_{H^s}\|D_{y'}v\|_{L^2}\\
&\qquad{}+\nu \big(1+(\theta\nu)^{-1/4}\big)\|w^0\|_{H^s}\|D_3 v\|_{L^2}\Big],
\end{align*}
the last inequality employing the fact that $\delta^{\alpha-1} \eta \leq K_0 \nu$, $\nu<\eta$, and $\theta\nu<1$. From this last inequality, we obtain that, for any $\varepsilon>0$,
\begin{align*}
& \left|\int_{\R^3_+} \wei{\A ( \nabla w+ \nabla \sB), \nabla v}\,dy\right| \\
&\leq C\Big[ \big( \eta+\eta(\theta\nu)^{1/4}+\nu(\theta\nu)^{-1/4}\big)\|w^0\|_{H^s}\|D_{y'}v\|_{L^2}\\
&\quad{}+\nu \big(1+(\theta\nu)^{-1/4}\big)\|w^0\|_{H^s}\|D_3 v\|_{L^2}\Big]\\
&\leq C\Big[ \eta\big( 1+(\theta\nu)^{1/4}+\frac{\nu}{\eta}(\theta\nu)^{-1/4}\big)\|w^0\|_{H^s}\|D_{y'}v\|_{L^2}\\
&\quad{}+\nu \big(1+(\theta\nu)^{-1/4}\big)\|w^0\|_{H^s}\|D_3 v\|_{L^2}\Big]\\
&\leq \eta \bigg(\frac{2C}{\varepsilon\Lambda}\big(1+\sqrt{\theta\nu}+\frac{\nu^2}{\eta^2}(\theta\nu)^{-1/2}\big)\|w^0\|_{H^s}^2+\varepsilon \cdot\frac{\Lambda}{2}\|D_{y'} v\|_{L^2}^2\bigg)\\
&\quad {}+\nu\bigg(\frac{2C}{\varepsilon\Lambda}\big(1+(\theta\nu)^{-1/2}\big)\|w^0\|_{H^s}^2+\varepsilon \cdot\frac{\Lambda}{2} \|D_3 v\|_{L^2}^2 \bigg)\\
&\leq \varepsilon\cdot \frac{\Lambda}{2}\left( \eta\| D_{y'}v\|_{L^2}^2+   \nu \| D_{3} v\|_{L^2}^2 \right)\\
&\quad{}+\frac{C}{\varepsilon}\left(\eta\big(1+\sqrt{\theta\nu}\big)+\nu\big(1+(\theta\nu)^{-1/2}\big)\right)\|w^0\|_{H^s}^2,
\end{align*}
where we have used that $\frac{\nu^2}{\eta}=\frac{\nu}{\eta}\nu\leq 1\cdot\nu$. The conclusion of the lemma follows from 
\[
\Lambda(\eta\|D_{y'}v\|_{L^2}^2+\nu\|D_3 v\|_{L^2}^2)=\int_{\R_+^3} \Lambda(\eta |D_{y'}v|^2+\nu |D_3 v|^2)\,dy
\]
and the first equation in statement (a) of Proposition \ref{prop:Z}.
\end{proof}

In our next lemma, we give the estimate for the fifth integral on the right hand side of \eqref{eq:Mas17}.
\begin{lemma}\label{le:Mas5}  Under the assumptions of \textup{Lemma \ref{le:Mas811}}, there exists a constant $C = C(L, s)>0$ such that
\begin{equation} \label{Est-le:Mas5}
\begin{split}
 \left|\int_{\R^3_+} \wei{([\bB v]\cdot\nabla)\sB,v}\,dy\right| 
& \leq  C\Big(\|w^0\|_{H^s}\|v\|_{L^2}^2 + \frac{\Lambda}{2}\varepsilon\eta \|D_{y'}v\|_{L^2}^2 \\
& \quad +\frac{2\theta \nu}{\Lambda\varepsilon \eta} \big(\|D_3 v\|_{L^2}^2+ \|v\|_{L^2}^2\big)\|w^0\|_{H^s}^2 \Big)
\end{split}
\end{equation}
holds for any $\varepsilon >0$ and for $\theta, \eta, \nu, \delta \in (0,1)$. In particular, if $\displaystyle \theta\leq \frac{\Lambda^2\varepsilon ^2\eta}{4\|w^0\|_{H^s}^2}$, then
\begin{equation} \label{rmk:c}
\begin{split}
& \left|\int_{\R^3_+} \wei{([\bB v]\cdot\nabla)\sB,v}\,dy\right|\\ 
& \leq C\left(\|w^0\|_{H^s}\|v\|_{L^2}^2+ \varepsilon \int_{\R^3_+} \wei{\A_0\bB^*\nabla v,\bB^*\nabla v}\,dy+\varepsilon \nu\|v\|_{L^2}^2\right),
\end{split}
\end{equation}
for any $\varepsilon>0$, and $\theta$, $\eta \in (0,1)$, $\nu \in (0, \eta)$, and for $\delta \in (0, \delta_0)$ such that $\delta^{\alpha-1} \eta \leq K_0 \nu$.
\end{lemma}
\begin{proof}
Let $j$, $k\in\{1,2\}$. We make the following claims:
\begin{align}
\|D_k\sB_j\|_{L^\infty}&\leq C\|w^0\|_{H^s},\label{eq:Mas5p1}\\
\|A(t,\cdot)\|_{L^\infty(\R^2)}&\leq C\|w^0\|_{H^s},\label{eq:Mas5p2}\\
\|D_k w_3(t,\cdot,0)\|_{L^\infty(\R^2)}&\leq C\|w^0\|_{H^s},\label{eq:Mas5p3}
\end{align}
where $C$ depends on $L$ and $s$, but is independent of $\delta$.

\smallskip
Equation \eqref{eq:Mas5p1} is a direct consequence of the first two lines of \eqref{BL.def}, \eqref{eq:w-def} (so $\sB_j$ does not involve $g$), the Sobolev embedding \eqref{eq:embed2d}, and the trace inequality \eqref{trace}. Equation \eqref{eq:Mas5p2} is similarly obtained from \eqref{eq:Aty}. As for \eqref{eq:Mas5p3}, we can expand $D_kw_3$ as follows:
\begin{multline*}
D_kw_3(t,y',0)=D_kw_3^0(t,y',0)+\delta^{\alpha-5/2}\bigg(D_kw_3^0(t,y',0)\\
{}+(\delta^{1/2}+\delta^{\alpha-2})\Big[D_kD_1g(y'/\delta)w_1^0(t,y',0)+D_kD_2g(y'/\delta)w_2^0(t,y',0)\Big]\\
{}+(\delta^{3/2}+\delta^{\alpha-1})\Big[D_1g(y'/\delta)D_k w_1^0(t,y',0)+D_2g(y'/\delta)D_k w_2^0(t,y',0)\Big]\bigg),
\end{multline*}
where $D_kD_1 g(y'/\delta)=\dfrac{\partial^2 g}{\partial y_k\partial y_1}\bigg|_{y'/\delta}$, and similarly for $D_kD_2 g(y'/\delta)$. As $\alpha>5/2$, we apply the triangle inequality, \eqref{eq:gL}, the Sobolev embedding \eqref{eq:embed2d}, and the trace inequality \eqref{trace} to yield \eqref{eq:Mas5p3}.

\smallskip
Once the preceding claims have been proven, the proof of Lemma 7 in \cite{PV} can be followed to obtain the conclusion \eqref{Est-le:Mas5}. The assertion \eqref{rmk:c} follows from the first inequality in Proposition \ref{prop:Z} (a) and \eqref{Est-le:Mas5}.
\end{proof}
The next lemma provides an estimate on the second to last line of \eqref{eq:Mas17}.
\begin{lemma}\label{le:newterms}  Let  $s, \alpha \in (5/2, \infty)$, and assume that \eqref{eq:gL} holds.  There is a constant $C = C(L, s)>0$ such that
\begin{align*}
& \left|\int_{\R_+^3}\wei{\partial_t \twd,v}\,dy\right|+\left|\int_{\R_+^3}\wei{(w^0\cdot\nabla)\twd,v}\,dy\right|+\left|\int_{\R_+^3}\wei{(\twd\cdot\nabla)w,v}\,dy\right|\\
& \leq C\left(\|w^0\|_{H^s}^2+\|F^0\|_{H^{s-1}}\right)\|v\|_{L^2}
\end{align*}
for all $\delta \in (0,1)$.
\end{lemma}
\begin{proof}
Using the first inequality in Lemma \ref{le:wprop} (c), we have
\[
\left|\int_{\R_+^3}\wei{\partial_t \twd,v}\,dy\right|\leq C\left(\|w^0\|_{H^s}^2+\|F^0\|_{H^{s-1}}\right) \|v\|_{L^2}.
\]
For the second term,
\begin{align*}
\left|\int_{\R_+^3}\wei{(w^0\cdot\nabla)\twd,v}\,dy\right|&\leq C\|(w^0\cdot\nabla) \twd\|_{L^2}\|v\|_{L^2}\\
&\leq C\|w^0\|_{H^s}^2\|v\|_{L^2}
\end{align*}
by the Banach algebra property \eqref{algebra-inq} and Lemma \ref{le:wprop} (e). The third term is estimated in a similar way.
\end{proof}
We finally provide an estimate for the last line of \eqref{eq:Mas17}.
\begin{lemma}\label{le:newterms2}  Let  $s, \alpha \in (5/2, \infty)$, and assume that \eqref{eq:gL} holds.  There is a constant $C= C(s, L)>0$ such that
\begin{align*}
& \delta^{\alpha-1}\left|\int_{\R_+^3}\wei{[(\bB_g w)\cdot\nabla]w,v}\,dy\right|+\left|\int_{\R_+^3}\wei{\nabla q,v}\,dy\right|\\
& \leq C\delta^{\alpha-1}\left(\|F^0\|+\|w^0\|_{H^s}^2\right)\|v\|_{L^2}
\end{align*}
for all $\delta \in (0, 1)$.
\end{lemma}
\begin{proof}
The estimate
\[
\left|\int_{\R_+^3}\wei{[(\bB_g w)\cdot\nabla]w,v}\,dy\right|\leq C\|w^0\|_{H^s}^2\|v\|_{L^2}
\]
can be deduced in the same way as the second term in Lemma \ref{le:newterms}. For the second integral, we first use \eqref{eq:bBexp} to write $\bI=\bB^*-\delta^{\alpha-1}\bB_g^*$, the divergence-free condition $\dive(\bB v)=0$, and $v\big|_{\{y_3=0\}}=0$, to find
\[
\int_{\R_+^3}\wei{\nabla q,v}\,dy=\int_{\R_+^3}\wei{(\bB^*-\delta^{\alpha-1}\bB_g^*)\nabla q,v}\,dy=-\delta^{\alpha-1}\int_{\R_+^3}\wei{\bB_g^*\nabla q,v}.
\]
On the other hand, from equation \eqref{eq:Euler} we obtain
\[
\nabla q=F^0-\partial_t w^0-(w^0\cdot\nabla)w^0.
\]
Taking $L^2$ norm on this equality, and applying the triangle inequality, Lemma \ref{le:wprop} (c), and the Banach algebra property \eqref{algebra-inq},
\begin{equation*}\label{eq:nablaq}
\|\nabla q\|_{L^2}\leq C\left( \|F^0\|_{L^2}+\|w^0\|_{H^s}^2\right),
\end{equation*}
whence
\[
\left|\int_{\R_+^3}\wei{\nabla q,v}\,dy\right|\leq C\delta^{\alpha-1}\left( \|F^0\|_{L^2}+\|w^0\|_{H^s}^2\right)\|v\|_{L^2}
\]
by \eqref{eq:bBg} and the boundedness of $Dg$ as in \eqref{eq:gL}. Combining the preceding estimates yields the conclusion of the lemma.
\end{proof}
Henceforth we fix   
\begin{equation}\label{eq:theta}
\theta=\frac{\Lambda^2\varepsilon^2\eta}{4\|w^0\|_{L^\infty((0,T^*),H^s)}^2+1}.
\end{equation}
Note that $\theta\nu<1$ for all $\eta$, $\nu\in (0,1)$ if $\varepsilon$ is sufficiently small, say, $0<\varepsilon\leq \varepsilon^*$ for some $\varepsilon^*>0$ depending on $\Lambda$ and $\|w^0\|_{L^\infty((0,T^*),H^s)}$. The value of $\varepsilon$ will be chosen below. Although the value of $\theta$ has been fixed, below we do not substitute its exact value, unless it is relevant to do so, for the sake of notational simplicity.

\smallskip
Using the previous lemmas, we can prove Proposition \ref{prop:G}.
\begin{proof}[Proof of Proposition \ref{prop:G}] Let $T^* \in (0, T)$ be as in Lemma \ref{Euler-existence}. We consider equation \eqref{eq:Mas17}. We first note that
\begin{equation}\label{eq:Fb}
\int_{\R^3_+} \wei{F-F^0,v}\,dy\leq \|F-F^0\|_{L^2} \|v\|_{L^2}.
\end{equation}
With $\theta$ as in \eqref{eq:theta} we see that  \eqref{rmk:c}  holds, so we can apply the previous lemmas, and \eqref{eq:Fb}, to equation \eqref{eq:Mas17} to conclude there is a constant $C^* = C^*(L, \Lambda, K_0, s)>0$ such that
\begin{equation}\label{eq:Gproof}
\begin{aligned}
& \frac{1}{2}\frac{d}{dt}\int_{\R^3_+} |v|^2\, dy + \int_{\R^3_+} \wei{\A_0 \bB^* \nabla v, \bB^* \nabla v}\, dy\\ 
& \leq C^*\Big[ \big((\theta\nu)^{1/4}+\delta^{\alpha-5/2}+\delta^{\alpha-1}\big)(\|w^0\|_{H^s}^2 + \|F^0\|_{H^{s-1}})\|v\|_{L^2} \\
&\quad {}+\|w^0\|_{H^s}\|v\|_{L^2}^2 +\varepsilon \int_{\R^3_+} \wei{\A_0\bB^*\nabla v,\bB^*\nabla v}\,dy+\varepsilon \nu\|v\|_{L^2}^2 \\ 
& \quad {}+\frac{1}{\varepsilon}\left(\eta(1+\sqrt{\theta\nu})+\nu(1+(\theta\nu)^{-1/2})\right)\|w^0\|_{H^s}^2\\
&\quad {}+\|F-F^0\|_{L^2} \|v\|_{L^2}\Big]
\end{aligned}
\end{equation}
holds for any positive constant $\varepsilon$, and for any  $\theta$, $\eta \in (0,1)$, $\nu \in (0, \eta)$, $\delta \in (0, \delta_0)$ satisfying $\delta^{\alpha-1} \eta \leq K_0 \eta$. Taking  $\varepsilon=\min\{\varepsilon^*,1/(2C^*)\}$, with $\varepsilon^*$ as in the remarks after \eqref{eq:theta}, we get
\[
C^* \varepsilon\int_{\R^3_+} \wei{\A_0\bB^*\nabla v,\bB^*\nabla v}\,dy\leq \frac{1}{2}\int_{\R^3_+} \wei{\A_0\bB^*\nabla v,\bB^*\nabla v}\,dy,
\]
the latter integral being nonnegative by the ellipticity of $\A_0$ (as in Proposition \ref{prop:Z} (a)), which means this term can be subtracted in the inequality \eqref{eq:Gproof}, and hence
\begin{equation}\label{eq:Mas18}
\begin{aligned}
& \frac{1}{2}\frac{d}{dt}\int_{\R^3_+} |v|^2\, dy+\int_{\R^3_+} \wei{\A_0 \bB^* \nabla v, \bB^* \nabla v}\, dy\\ 
& \leq C\Big[\Big(\|w^0\|_{H^s}+\nu\big)\|v\|_{L^2}^2\\ 
& \quad {}+\Big(\big((\theta\nu)^{1/4}+\delta^{\alpha-5/2}+\delta^{\alpha-1}\big)(\|w^0\|_{H^s}^2 + \|F^0\|_{H^{s-1}})\\
&\qquad {} +\|F-F^0\|_{L^2}\Big) \|v\|_{L^2}\\
& \quad{}+\left(\eta(1+\sqrt{\theta\nu})+\nu(1+(\theta\nu)^{-1/2})\right)\|w^0\|_{H^s}^2\Big].
\end{aligned}
\end{equation}
Let
\begin{align*}
f_0(t)&=C\big(\|w^0\|_{H^s}+\nu\big),\\
f_1(t)&= C\Big(\big((\theta\nu)^{1/4}+\delta^{\alpha-5/2}\big)(\|w^0\|_{H^s}^2 + \|F^0\|_{H^{s-1}})\\
&\qquad {} +\|F-F^0\|_{L^2}\Big),\\
f_2(t)&=C\left(\eta+\sqrt{\nu/\eta} \|w^0\|_{L^\infty((0,T^*),H^s)} \right)\|w^0\|_{H^s}^2.
\end{align*}
Our choice of $f_1(t)$ arises from the fact that $\delta^{\alpha-1}<\delta^{\alpha-5/2}$ (since $\delta<1$), while our choice of $f_2(t)$ is explained by noting that, under our assumptions, $\theta\nu<1$, $\nu<\eta$, and $\nu(\theta\nu)^{-1/2} \leq C\|w^0\|_{L^\infty((0,T^*),H^s)}\sqrt{\nu/\eta}$, and consequently
\[
C\left(\eta(1+\sqrt{\theta\nu})+\nu(1+(\theta\nu)^{-1/2})\right) \|w^0\|_{H^s}^2 \leq f_2(t).
\]
Clearly \eqref{eq:ineq} holds with our choices of $f_0$, $f_1$, and $f_2$,  and $f_0$, $f_1$, $f_2\in L^1(0,T^*)$, since $\|w^0(t,\cdot)\|_{H^s}\in L^\infty(0,T^*)$, $\|F^0(t,\cdot)\|_{H^{s-1}}\in L^1(0,T^*)$, and $\|F(t,\cdot)-F^0(t,\cdot)\|_{L^2}\in L^1(0,T^*)$  as assumed in \eqref{ini-cond}. From $w^0\in L^\infty((0,T^*),H^s(\R_+^3)^3)$ and $0<\nu<1$, we see that
\[
\int_0^{T^*} f_0(t)\,dt\leq C.
\]
Similarly, we find that
\[
\int_0^{T^*} f_1(t)\,dt\leq C\big((\theta\nu)^{1/4}+\delta^{\alpha-5/2}+\beta(\eta,\nu)\big)\leq C(\beta(\eta,\nu)+\delta^{\alpha-5/2}),
\]
where we have used the second inequality in \eqref{ini-cond} and the fact that $(\theta\nu)^{1/4}=C(\eta\nu)^{1/4}\leq C(\nu/\eta)^{1/4}\leq C\beta(\eta,\nu)$, since $\nu<\eta<1$. Finally,
\[
\int_0^{T^*} f_2(t)\,dt\leq C\beta^2\leq C(\beta+\delta^{\alpha-5/2})^2,
\]
so the inequalities in \eqref{eq:abcboundsbeta} are all satisfied.
\end{proof}
\begin{rmkn}\label{rmk:abc}
We note that the constant $C$ involved in the definition of the functions $f_0$, $f_1$, and $f_2$ depends on $L$ \textup{(}as in \eqref{eq:gL}\textup{)}, $\|w^0\|_{L^\infty((0,T^*),H^s)}$, $\|F^0\|_{L^1((0,T^*),H^s(\R_+^3))}$, $T^*, K_0$, and $\Lambda$. We also note that $\varepsilon$ (and thus $\theta$) has been fixed, $(\theta\nu)^{1/4}\leq C\beta(\eta,\nu)$, and $\|w^0\|_{L^\infty((0,T^*),H^s)}\leq C$, which will be used below.
\end{rmkn}
\section{Proof of Theorem \ref{thm:main}} \label{proof-thrm}
Before we prove our main theorem, we state the following Gr\"{o}nwall type lemma, which is proved in \cite[Lemma 1.1]{Mas}.
\begin{lemma} \label{le:Gronwall}
Let $\gamma>0$, $T \in (0, \infty)$, and let $f_0$, $f_1$, $f_2\in L^1(0,T)$ be nonnegative functions satisfying
\begin{equation}\label{eq:abcbounds}
\int_0^T f_0(t)\,dt\leq C ,\quad \int_0^T f_1(t)\,dt\leq C\gamma ,\quad \int_0^T f_2(t)\,dt\leq C\gamma^2
\end{equation}
for some constant $C>0$. Then, there exists a constant $M = M(C)>0$ such that for every nonnegative function $f$ satisfying
\begin{equation}\label{eq:dtf}
\partial_t (f^2)\leq f_0(t)f^2+f_1(t)f+f_2(t) \quad t \in (0, T) \quad \text{and}\quad f(0)\leq C\gamma,
\end{equation}
it holds that
\[
|f(t)|\leq M\gamma, \quad \text{for all }t \in [0,T].
\]
\end{lemma}
We now have all the ingredients for the proof of Theorem \ref{thm:main}.
\begin{proof}[Proof of Theorem \ref{thm:main}]  Let $f(t)=\|v(t,\cdot)\|_{L^2}$, $t\in (0,T^*)$. Our goal is to apply Lemma \ref{le:Gronwall} to $f$ with $C$, $f_0$, $f_1$, and $f_2$ as provided by Proposition \ref{prop:G}, and $\gamma=\beta(\eta,\nu)+\delta^{\alpha-5/2}$. Under these choices, equation \eqref{eq:abcbounds} is a consequence of \eqref{eq:abcboundsbeta}, while the first inequality in \eqref{eq:dtf} is obtained from \eqref{eq:ineq} by dropping the second term on the left hand side, which is nonnegative due to the ellipticity of $\A_0$ (see Proposition \ref{prop:Z} (a)). It remains to verify $f(0)\leq C(\beta(\eta,\nu)+\delta^{\alpha-5/2})$. First, using \eqref{eq:sB}, Lemma \ref{le:bounds1} (a) and (b), Lemma \ref{le:wprop} (d), and Remark \ref{rmk:abc}, we have 
\begin{equation}\label{eq:sBinftybound2}
\|\sB(0,\cdot)\|_{L^2}\leq C(\theta\nu)^{1/4}\|w^0(0,\cdot)\|_{H^s} \leq C\beta(\eta,\nu).
\end{equation}
Then, recalling that $v=u-w-\sB$, using \eqref{ini-cond}, \eqref{eq:w-def}, and \eqref{eq:sBinftybound2}, we obtain
\begin{align*}
f(0)& =\|u(0,\cdot)-w(0,\cdot)-\sB(0,\cdot)\|_{L^2}\\
&\leq \|U (\cdot)-(w^0(0,\cdot)-\delta^{\alpha-5/2}\twd(0,\cdot ))\|_{L^2}+\|\sB(0,\cdot)\|_{L^2}\\
&\leq \|U (\cdot)-W^0(\cdot)\|_{L^2}+\delta^{\alpha-5/2}\|\twd(0,\cdot)\|_{L^2}+\|\sB(0,\cdot)\|_{L^2}\\
&\leq C(\beta(\eta,\nu)+\delta^{\alpha-5/2}),
\end{align*}
where we have used that $\|\twd(0,\cdot)\|_{L^2}\leq C\|w^0\|_{L^2}\leq C$ (see Remark \ref{rmk:abc}). 

\smallskip
Having verified all the assumptions of Lemma \ref{le:Gronwall}, there exists a constant $M>0$, depending on $L$, $\|w^0\|_{L^{\infty}((0,T^*), H^s)}$, $\Lambda, K_0$, and $T^*$, such that
\begin{equation}\label{eq:vbound}
\|v\|_{L^\infty((0,T^*),L^2(\R^3_+))}\leq M(\beta(\eta,\nu)+\delta^{\alpha-5/2}).
\end{equation}

\smallskip
Now, we prove that \eqref{eq:uwbound1} holds. From Proposition \ref{prop:sBbounds} (e) and Remark \ref{rmk:abc}, we have
\begin{align}
\nonumber \|\sB\|_{L^\infty((0,T^*),L^2(\R^3_+))} & \leq C\|\bB\sB\|_{L^\infty((0,T^*),L^2(\R^3_+))}\\
\nonumber &\leq C(\theta\nu)^{1/4}\|w^0\|_{L^\infty((0,T^*),H^s(\R^3_+))}\\
&\leq C\beta(\eta,\nu).\label{eq:sBinftybound}
\end{align}
Then, as $u=v+w+\sB$, it follows from \eqref{eq:vbound}, \eqref{eq:sBinftybound}, and the triangle inequality that
\[
\|u-w\|_{L^\infty((0,T^*),L^2(\R^3_+))}\leq C(\beta+\delta^{\alpha-5/2}),
\] 
where $C$ depends on $T^*$, $\|w^0\|_{L^\infty((0,T^*), H^s(\R_+^3))}$, $\Lambda, K_0$, and $L$. Applying the triangle inequality again, together with \eqref{eq:w-def}, and $\|\twd\|_{L^2}\leq C\|w^0\|_{L^2}$,
\[
\|u-w^0\|\leq \|u-w\|+\|w-w^0\|\leq C(\beta+\delta^{\alpha-5/2}),
\]
all norms taken on ${L^\infty((0,T^*),L^2(\R^3_+))}$. From this last estimate  we conclude that equation \eqref{eq:uwbound1} holds.

\smallskip
It remains to prove that \eqref{eq:uwbound2} is satisfied. Integrating \eqref{eq:ineq} with respect to the time variable on $(0, T^*)$ and using  \eqref{eq:abcboundsbeta} and \eqref{eq:vbound}, we find
\[
\int_0^{T^*}\int_{\R^3_+} \wei{\A\nabla v(t, y),\nabla v(t, y)}\, dydt \leq C(\beta+\delta^{\alpha-5/2})^2.
\]
From Proposition \ref{prop:Z} (a), it follows that
\begin{equation}  \label{eq:thmproof1}
\int_0^{T^*}\int_{\R^3_+} \left(\eta|D_{y'}v|^2+\nu|D_3v|^2\right)\, dydt \leq C(\beta+\delta^{\alpha-5/2})^2.
\end{equation}
Recalling that $0<\nu<\eta<1$, $\|w\|_{L^\infty((0,T^*),H^1)}\leq C\|w^0\|_{L^\infty((0,T^*),H^s)}\leq C$ (by Lemma \ref{le:wprop} (e)), a computation yields
\begin{align}
\int_0^{T*}\int_{\R_+^3}\left(\eta|D_{y'}w|^2+\nu|D_3w|^2\right)\,dydt&\leq C(\eta+\nu)\nonumber\\
&\leq C(\eta+\sqrt{\nu/\eta}) \nonumber\\
&\leq C\beta^2.\label{eq:thmproof2}
\end{align}
As $\|w^0\|_{L^\infty((0,T^*),H^s)} \leq C$ and $\theta\nu<1$, we obtain from Proposition \ref{prop:sBbounds} (b) that  
\[
\eta\|D_{y'}\sB(t,\cdot)\|_{L^2}^2\leq C\eta,
\]
and, using again that $\theta\nu=C\eta\nu$ and Proposition \ref{prop:sBbounds} (c),
\[
\nu\|D_3\sB(t,\cdot)\|_{L^2}^2\leq C\nu(\theta\nu)^{-1/2}=C\sqrt{\nu/\eta}.
\]
Combining the two foregoing estimates, and integrating on $(0,T^*)$, we find
\begin{equation}\label{eq:thmproof3}
\int_0^{T*}\int_{\R^3_+} \left(\eta|D_{y'}\sB|^2+\nu|D_3\sB|^2\right)\ dydt\leq C(\eta+\sqrt{\nu/\eta})\leq C\beta^2.
\end{equation}
Using the second inequality in Proposition \ref{prop:Z} (a), the triangle inequality, and the estimates \eqref{eq:thmproof1}, \eqref{eq:thmproof2}, \eqref{eq:thmproof3},
\begin{align*}
&\int_0^{T^*}\int_{\R_+^3} \wei{\A(t,y)\nabla u(t,y),\nabla u(t,y)}\,dydt\\
&\leq C\int_0^{T^*}\int_{\R_+^3} \left[\eta |D_{y'} u|^2+\nu|D_3 u|^2\right]\,dydt\\
&\leq C\int_0^{T^*}\int_{\R_+^3} \Big[\eta (|D_{y'} v|^2+|D_{y'} w|^2+|D_{y'} \sB|^2)\\
&\qquad {}+\nu(|D_3 v|^2+|D_3 w|^2+|D_3 \sB|^2)\Big]\,dydt\\
&\leq C(\beta(\eta,\nu)+\delta^{\alpha-5/2})^2.
\end{align*}
By the change of variables $y=\Phi_0(x)$,
\begin{equation*}
\int_0^{T^*}\int_{\Omega} \wei{\A_0\nabla \bar{u}(t, x),\nabla \bar{u}(t, x)}\, dx dt\leq C (\beta+\delta^{\alpha-5/2})^2,
\end{equation*}
where $C$ depends on $T^*$, $\Lambda$, $L$, $K_0$, $s$, and $\|w^0\|_{L^\infty((0,T^*),H^s(\R_+^3))}$. Then, the assertion \eqref{eq:uwbound2} follows by the preceding estimate and the first assumption in \eqref{ellipticity-cond}. This concludes the proof.
\end{proof}

\end{document}